\documentclass{elsarticle}
				
\usepackage[a4paper]{geometry} 
\usepackage[utf8]{inputenc}
\usepackage{amsmath}
\usepackage{amssymb}
\usepackage{bm} 
\usepackage{graphicx}
\usepackage{epstopdf} 
\usepackage{multirow}
\usepackage{tocbasic}
\usepackage{caption}
\usepackage{subcaption}
\usepackage{placeins}
\usepackage[shortcuts]{extdash}
\usepackage{siunitx}
\usepackage{mdframed}
\usepackage{url}
\usepackage{framed}
\usepackage{lscape} 
\usepackage{booktabs} 
\usepackage{todonotes} 
\usepackage{empheq}

\graphicspath{{figures/}}
\usepackage{hyperref} 

\newcommand{\vt}[1]{\bm{#1}}
\newcommand{\rd}{\textrm{d}}

\newcommand{\dd}[2]{\frac{\partial #1}{\partial #2}}

\newcommand{\hl}[1]{\textcolor{black}{#1}}
\newcommand{\hll}[1]{\textcolor{black}{#1}}

\usepackage{helvet}

\usepackage{lmodern}

\begin{document}

\title{Constraint-consistent Runge-Kutta methods for one-dimensional incompressible multiphase flow}    
\author[add1]{B. Sanderse}
\ead{B.Sanderse@cwi.nl, corresponding author}
\author[add2]{A.E.P. Veldman}
\address[add1]{Centrum Wiskunde \& Informatica (CWI), Amsterdam, The Netherlands}
\address[add2]{Bernoulli Institute for Mathematics, Computer Science and Artificial Intelligence, University of Groningen, Groningen, The Netherlands}

\begin{abstract}

New time integration methods are proposed for simulating incompressible multiphase flow in pipelines described by the one-dimensional two-fluid model. The methodology is based on `half-explicit' Runge-Kutta methods, being explicit for the mass and momentum equations and implicit for the volume constraint. These half-explicit methods are \textit{constraint-consistent}, i.e., they satisfy the hidden constraints of the two-fluid model, namely the volumetric flow (incompressibility) constraint and the Poisson equation for the pressure. A novel analysis shows that these hidden constraints are present in the continuous, semi-discrete, and fully discrete equations.  

Next to constraint-consistency, the new methods are \textit{conservative}: the original mass and momentum equations are solved, and the proper shock conditions are satisfied; \textit{efficient}: the implicit constraint is rewritten into a pressure Poisson equation, and the time step for the explicit part is restricted by a CFL condition based on the convective wave speeds; and \textit{accurate}: achieving high order temporal accuracy for all solution components (masses, velocities, and pressure). High-order accuracy is obtained by constructing a new third-order Runge-Kutta method that satisfies the additional order conditions arising from the presence of the constraint in combination with time-dependent boundary conditions.

Two test cases (Kelvin-Helmholtz instabilities in a pipeline and liquid sloshing in a cylindrical tank) show that for time-independent boundary conditions the half-explicit formulation with a classic fourth-order Runge-Kutta method accurately integrates the two-fluid model equations in time while preserving all constraints. A third test case (ramp-up of gas production in a multiphase pipeline) shows that our new third-order method is preferred for cases featuring time-dependent boundary conditions.

\end{abstract}



\maketitle

\section{Introduction}

\subsection{Background}
The incompressible two-fluid model for one-dimensional multiphase flow is an important model to study, for example, the behaviour of oil and gas in long pipelines. A main research area is the development of methods that accurately solve the two-fluid model in order to predict the transition from stratified flow to slug flow, and the subsequent propagation of the generated slugs (so-called slug capturing methods). In this paper we progress towards this goal by developing efficient time integration methods for the incompressible two-fluid model in its basic form: \hl{a four-equation model (describing conservation of mass and momentum per phase), supplemented with a volume constraint (describing that the phases together exactly fill the pipeline)}. An important aspect in the construction of time integration methods lies in the correct treatment of this volume constraint and associated derived constraints, such as the divergence-free constraint on the mixture velocity field.

Historically, the development of efficient time integration methods for the two-fluid model hinges on the use of the pressure equation. For the \textit{compressible} two-fluid model a hyperbolic evolution equation for the pressure can be derived by a single differentiation of the volume constraint (see e.g.\ \cite{Bendiksen1991, Moe1993}). The resulting system can be solved for example with implicit methods in order to circumvent the CFL condition associated to the acoustic waves \cite{Evje2005,Sanderse2017}. However, in case the flow is incompressible, the character of the pressure equation changes from hyperbolic to elliptic, and dedicated incompressible solvers are more efficient than compressible solvers.

Different approaches to deal with this incompressible, one-dimensional, multiphase flow problem have been proposed, which will be shortly summarized here. 

A first approach is to eliminate the pressure from the four-equation system and to rewrite this system into a two-equation system. This leads to the `no-pressure wave' model or the `fixed-flux' model suggested by \cite{Watson1989}, and used for example in \cite{Akselsen2016,Holmas2008,LopezdeBertodano2017,Omgba-Essama2004}. A similar two equation model is the `reversed density' model developed by Keyfitz et al.\ \cite{Keyfitz2003} and employed for example in \cite{Wangensteen2010}. In these models, the pressure is generally computed as a post-processing step. A general problem with these reduced equation systems is that they are only valid in case of smooth solutions. In the presence of shocks the wrong jump conditions are obtained \cite{Akselsen2016}. Furthermore, the fixed-flux assumption often limits these studies to stationary boundary conditions. 

A second approach is to keep the pressure in the formulation and to use a pressure-correction method. A pressure equation is then typically obtained by substituting the momentum equations in the combined mass conservation equation, while applying the volume constraint equation. This approach is taken by Liao et al.\ \cite{Liao2008}, who solve the momentum equations with an old pressure, solve for the new pressure, and then update the velocity. A slightly different approach is taken in \cite{Issa2003,LopezdeBertodano2017,Wang2015}, in which the pressure equation is derived from a combination of momentum and continuity equations, without using the constraint equation. A related approach is the all-speed method in the RELAP-7 code \cite{Berry2014}, based on the PCICE (pressure corrected implicit continuous-fluid eulerian) algorithm. The temporal accuracy of these approaches is mostly limited to first order, or not reported. Note that most of these methods are reminiscent of single phase incompressible (2D or 3D) Navier-Stokes solution algorithms such as SIMPLE (\cite{Patankar1980}) or PISO \cite{Issa1986}, which were developed to handle the divergence-free constraint. 

The aim of this paper is to resolve the issues of these current approaches by developing high order \textit{constraint-consistent} time integration methods for the incompressible two-fluid model equations in conservation form, including a generic (non-stationary) boundary condition treatment.

\subsection{Approach and outline}

In this paper we study the incompressible two-fluid model equations in conservative form. The constraint in the model has implications on all discretization aspects and we have constructed the roadmap of figure \ref{fig:DAE_ODE} to make this clear. First we show via a characteristic analysis that the constraint in the continuous model equations introduces two additional `hidden' constraints (section \ref{sec:governing_equations}, first row of figure \ref{fig:DAE_ODE}). Second, the constraint has important consequences for the spatial discretization of the equations, in particular for the boundary conditions, which are derived in section \ref{sec:spatial_discretization_bc}. Third, in the resulting semi-discrete equations the presence of the constraint leads to a differential-algebraic equation (DAE) system that features the same two hidden constraints (section \ref{sec:temporal_properties}, second row in figure \ref{fig:DAE_ODE}). These semi-discrete equations are discretized in time with a `half-explicit' Runge-Kutta method and recombined into a form that highlights the hidden constraints in the equations (section 5, third row in figure \ref{fig:DAE_ODE}). The DAE classification is used to perform the accuracy analysis, and a new third-order Runge-Kutta method will be designed such to avoid order reduction. The resulting method is a new high-order time integration method that is consistent with the constraints derived on the continuous and semi-discrete level.  In summary, our novel approach follows the red dashed line in figure \ref{fig:DAE_ODE} and is consistent with the mantra: \textit{discretize first, substitute next} \cite{Veldman1990}. 

The main advantages of this approach are that (i) the method is consistent with the constraints of the model, (ii) high-order accuracy of all solution components, including the pressure, is achieved, (iii) the original conservation equations are solved, and (iv) the approach only requires implicit treatment of the pressure equation. Furthermore, the approach has the potential to be applied to other constraint systems, such as drift-flux models and three-fluid models, and multi-dimensional applications. 

\begingroup
\begin{figure}[hbtp]
\fontfamily{lmss} 
\fontsize{9pt}{10pt}\selectfont
\centering 
\def\svgwidth{1.0 \textwidth}
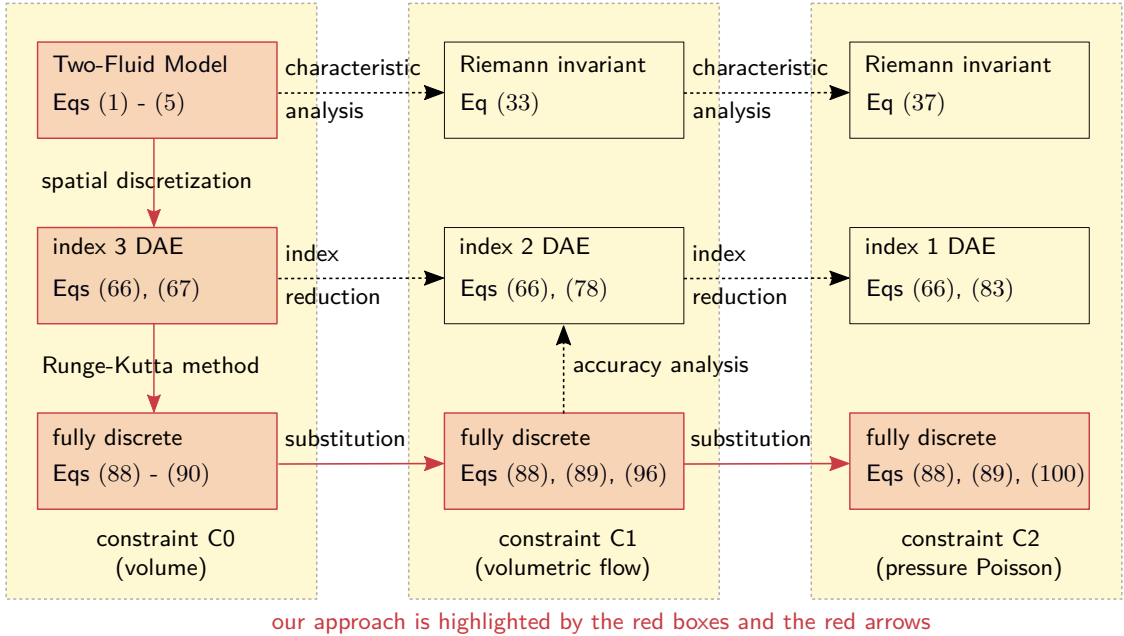 
\caption{Our approach to form new time integration methods for the incompressible two-fluid model follows the red boxes: \textit{discretize first, substitute next}.}
\label{fig:DAE_ODE}
\end{figure}
\endgroup


\section{Constraint analysis of differential equations}\label{sec:governing_equations}

\subsection{Governing equations incompressible flow}
The incompressible two-fluid model can be derived by considering the stratified flow of liquid and gas in a pipeline (for a recent discussion of the two-fluid model, see for example \cite{LopezdeBertodano2017}). \hl{For a sketch of the geometry, see figure \ref{fig:stratified_flow}.}
\begingroup
\begin{figure}[hbtp]
\fontfamily{lmss}
\fontsize{10pt}{12pt}\selectfont
\centering 
\def\svgwidth{0.8 \textwidth}
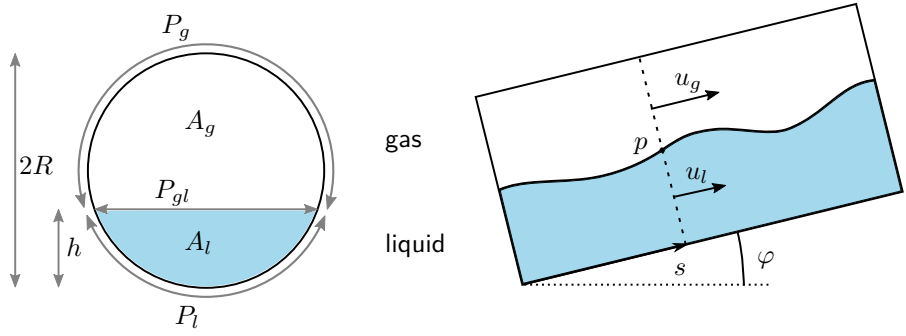 
\caption{\hl{Stratified flow in a pipeline. Left: cross-sectional view; right: side view.}}
\label{fig:stratified_flow}
\end{figure}
\endgroup
The main assumptions that we make are that the flow is one-dimensional, stratified, incompressible, and isothermal. Transverse pressure variation is introduced via level gradient terms. Surface tension is neglected. This leads to the following \hl{four-equation} model:
\begin{align}
\frac{\partial}{\partial t} \left( \rho_{g} A_{g} \right) + \frac{\partial }{\partial s} \left( \rho_{g} u_{g} A_{g} \right) &= 0, \label{eqn:conservation_mass_gas} \\
\frac{\partial}{\partial t} \left( \rho_{l} A_{l} \right) + \frac{\partial }{\partial s} \left( \rho_{l} u_{l} A_{l} \right) &= 0,\label{eqn:conservation_mass}\\
\frac{\partial }{\partial t} \left( \rho_{g} u_{g} A_{g} \right) + \frac{\partial }{\partial s} (\rho_{g} u_{g}^2 A_{g}) &= -\frac{\partial p }{\partial s} A_{g} + E_{g} + \underbrace{ - \tau_{gl} P_{gl} - \tau_{g} P_{g} - \rho_{g} A_{g} g_s +  F_{\text{body}} {A_g}}_{S_{g}},\\
\frac{\partial }{\partial t} \left( \rho_{l} u_{l} A_{l} \right) + \frac{\partial }{\partial s} (\rho_{l} u_{l}^2 A_{l}) &= -\frac{\partial p }{\partial s} A_{l} + E_{l} + \underbrace{\tau_{gl} P_{gl} - \tau_{l} P_{l} - \rho_{l} A_{l} g_s +  F_{\text{body}} {A_l}}_{S_{l}},\label{eqn:conservation_momentum}
\end{align}
supplemented with the \textit{volume constraint equation}, which we indicate by C0:
\begin{equation}\label{eqn:volume_constraint}
\text{\underline{C0}:} \qquad A_g + A_l = A.
\end{equation}
In these equations the subscript denotes either gas ($g$) or liquid ($l$). The model features four evolution equations, one constraint equation, and five unknowns ($A_g$, $A_l$, $u_{g}$, $u_l$, $p$), which are a function of the independent variables $s$ (coordinate along the pipeline axis) and $t$ (time). $\rho$ denotes the density (assumed constant), $A$ the cross-sectional area of the pipe, $A_{g}$ and $A_{l}$ (also referred to as the hold-ups) the cross-sectional areas occupied by the gas or liquid, $u$ the phase velocity, $p$ the pressure at the interface, $\tau$ the shear stress (with the wall or at the interface), $g$ the gravitational constant, $\varphi$ the local inclination of the pipeline with respect to the horizontal, $g_n=g \cos \varphi$ and $g_s = g \sin \varphi$. \hl{The wetted and interfacial perimeters $P_{g}$, $P_{l}$ and $P_{gl}$ can be expressed in terms of the hold-ups $A_{g}$ or $A_{l}$ via a non-linear algebraic expression, see \ref{sec:geometry}. The same is true for the interface height $h$ (measured from the bottom of the pipe), which appears in the expression for the level gradient terms $E$ \cite{Sanderse2017}:
\begin{align}\label{eqn:compressible_HG}
E_{g} &= \frac{\partial K_{g}}{\partial s}, \quad K_{g} = \rho_{g} g_n \left[ (R-h) A_{g} + \frac{1}{12} P_{gl}^3\right], \\
E_{l} &= \frac{\partial K_{l}}{\partial s}, \quad K_{l} = \rho_{l} g_n \left[ (R-h) A_{l} - \frac{1}{12} P_{gl}^3\right].
\end{align}}\hl{In these expressions $R$ is the pipe radius}. For incompressible flow these expressions simplify to $E_{g}=-\rho_{g} g_{n} A_{g} \frac{\partial h}{\partial s}$, $E_{l}=-\rho_{l} g_{n} A_{l} \frac{\partial h}{\partial s}$, but we stick to the form displayed in equation \eqref{eqn:compressible_HG} because this form is conservative. The two-fluid model can also be employed for channel flow instead of pipe flow with minor modifications.

The body force $F_{\text{body}}$ in the gas and liquid momentum equations is for example a driving pressure force for the simulations that involve periodic boundary conditions, or a source term to force an analytical solution (see \ref{sec:MMS}). The friction models for $\tau_{g}$, $\tau_{l}$ and $\tau_{gl}$ are described in \ref{sec:friction_models}. The source terms $S_{g}$ and $S_{l}$ do not contain spatial or temporal derivatives.

\hll{The two-fluid model is supplemented with initial and boundary conditions to form an initial boundary value problem. As the two-fluid model equations are conditionally hyperbolic \cite{Barnea1994,Lyczkowski1978,Stewart1984}, not all initial and boundary conditions guarantee that the equations remain hyperbolic. In this paper we will avoid non-hyperbolic conditions and consequently guarantee well-posedness, although we should note that for non-hyperbolic equations results have also been obtained \cite{Keyfitz2003}.}

\subsection{Two-equation model}\label{sec:alternative_forms}
Several discretization methods for solving the two-fluid model equations \eqref{eqn:conservation_mass_gas}-\eqref{eqn:volume_constraint} are based on an alternative form of the model. For the compressible model, a common approach is to expand the density in the mass conservation equations in terms of the pressure, and combine the resulting equations with the volume constraint to derive a pressure equation \cite{Bendiksen1991, Evje2005, Evje2007}. However, in case the fluid is incompressible, such a pressure equation becomes singular and cannot be used. 
A common approach is then to reduce the four-equation system to a two-equation system \cite{Akselsen2016,Holmas2008,Omgba-Essama2004}. The basis lies in combining the mass equations and the volume constraint to derive the \textit{volumetric flow constraint} (C1):
\begin{equation}\label{eqn:mixture_velocity_constraint1}
\text{\underline{C1}:} \qquad \frac{\partial }{\partial s} \left( A_{g} u_{g} + A_{l} u_{l} \right) = 0.
\end{equation}
It can be integrated in space to give, at any point in the pipeline, 
\begin{equation}\label{eqn:mixture_velocity_constraint2}
A u_{\text{mix}} := A_{g} u_{g} + A_{l} u_{l} =  V(t),
\end{equation}
where $V(t)$ is the volumetric flow rate. The next step is to rewrite the momentum and mass equations into an equation for the difference in momenta:
\begin{equation}\label{eqn:combined_momentum}
\frac{\partial }{\partial t} \left( \rho_{l} u_{l} - \rho_{g} u_{g} \right) + \frac{\partial }{\partial s} \left( \frac{1}{2} \left(\rho_{l} u_{l}^2 - \rho_{g} u_{g}^2 \right) \right) 
= - (\rho_{l} - \rho_{g}) g_{n} \frac{\partial h}{\partial s} + \frac{S_{l}}{A_{l}}- \frac{S_{g}}{A_{g}}. 
\end{equation}
At the same time, the mass equations can be combined to give
\begin{equation}\label{eqn:combined_mass}
\frac{\partial}{\partial t} \left( \rho_{g} A_{g} + \rho_{l} A_{l}\right) + \frac{\partial }{\partial s} \left( \rho_{g} u_{g} A_{g} +\rho_{l} u_{l} A_{l}  \right) = 0.
\end{equation}
The `conservative' (pressure-free) combined momentum equation \eqref{eqn:combined_momentum} and the conservative combined mass equation \eqref{eqn:combined_mass} form a two-equation system, which can be numerically solved when supplemented with the volume constraint \eqref{eqn:volume_constraint} and volumetric flow constraint \eqref{eqn:mixture_velocity_constraint2}. 
However, this reduced pressure-free model is only equivalent to the original two-fluid model for sufficiently smooth solutions. In case shock waves appear in the solution, the pressure-free model is not equivalent, and will generally exhibit different weak solutions \cite{Akselsen2016a,LeVeque2002}. An example of shock waves in the two-fluid model are roll waves \cite{Watson1989}, for which it has been shown that the conservative and non-conservative equation systems indeed lead to different solutions \cite{Akselsen2016a,Sanderse2018b}. Consequently, we will \textit{not} pursue the two-equation system as a starting point for the development of time integration methods, but instead use the original conservation equations \eqref{eqn:conservation_mass_gas}-\eqref{eqn:conservation_momentum} \footnote{We note that the momentum equations in the four-equation model are not fully conservative due to the presence of the hold-up fractions in the pressure gradient terms. \hll{This is a topic of study in itself, and we refer the interested reader to \cite{Castro2008,Toumi1996}.}}.

\subsection{Constraint equations from characteristic analysis}
The governing equations of the two-fluid model describe the change in mass and momentum of each phase, commonly expressed by the vector of conservative variables
\begin{equation}
\vt{U} = \begin{pmatrix}
\rho_{g} A_{g} \\
\rho_{l} A_{l} \\
\rho_{g}  u_{g} A_{g} \\
\rho_{l} u_{l} A_{l}
\end{pmatrix}.
\end{equation}
As mentioned, the equations cannot be written in full conservative form, but fortunately they can be written in quasi-linear form in terms of the primitive variables $\vt{W}$,
\begin{equation}\label{quasi-linear}
\vt{A} (\vt{W}) \frac{\partial \vt{W}}{\partial t}  + \vt{B}( \vt{W}) \frac{\partial \vt{W}}{\partial s} + \vt{S}(\vt{W}) = \vt{0},
\end{equation}
where we choose  $\vt{W} = (A_{l}, u_{l}, u_{g}, p)^T$.
The dependence of $\vt{A}$, $\vt{B}$ and $\vt{S}$ on $\vt{W}$ will not be explicitly indicated in the sequel. The matrices $\vt{A}$ and $\vt{B}$ and source term $\vt{S}$ can be simplified, by multiplying the mass equations by their respective velocities and adding them to the momentum equations, and by dividing by the respective cross-sectional areas:
\begin{equation}
\vt{A} = 
\begin{pmatrix}
1 & 0 & 0 & 0 \\
-1 & 0 & 0 & 0 \\
0 & \rho_{l} & 0 & 0 \\
0 & 0 & \rho_{g} & 0 \\
\end{pmatrix},\qquad
\vt{B} = 
\begin{pmatrix}
u_{l} & A_{l} & 0 & 0\\
- u_{g} & 0 & A_{g} &  0\\
- \frac{1}{A_{l}} \frac{\partial K_{l}}{\partial A_{l}}  &  \rho_{l} u_{l}  & 0 & 1 \\
+ \frac{1}{A_{g}} \frac{\partial K_{g}}{\partial A_{g}}   & 0 & \rho_{g} u_{g} & 1
\end{pmatrix},
\qquad
\vt{S} = \begin{pmatrix}
0 \\
0 \\
- S_{l}/A_{l} \\ 
-S_{g}/A_{g} 
\end{pmatrix}.
\end{equation}
In this derivation the constraint has been substituted, in the form $\dd{A_{g}}{t} = -\dd{A_{l}}{t}$ and $\dd{A_{g}}{s} = -\dd{A_{l}}{s}$. The derivatives of the level gradient terms read
\begin{align}
\frac{\partial K_{l}}{\partial A_{l}}  &= -\rho_{l} g_{n} H_{l}, \qquad \frac{\partial K_{g}}{\partial A_{g}}  = \rho_{g} g_{n} H_{g},
\end{align}
where $H_{l}/A_{l} = H_{g}/A_{g} = \dd{h}{A_{l}}$ (similar to the derivation in \cite{Barnea1994}). 

The eigenvalues are found from the generalized eigenvalue problem 
\begin{equation}
\text{det} ( \vt{B} - \lambda \vt{A}) = 0.
\end{equation}
Since $\vt{A}$ is rank-deficient ($\text{rank} (\vt{A}) = 3$), the number of eigenvalues is lower than the dimension of the matrix \cite{Golub2012}. Two of the roots are given by 
\begin{equation}\label{eqn:roots}
\lambda_{1,2} = \frac{ (\rho u)^{*} \pm \xi}{ \rho^{*}},
\end{equation}
where, in the notation from Akselsen \cite{Akselsen2016a}, \hll{we define the averaging operator $(.)^*$ as:}
\begin{equation}
(.)^{*} := \frac{ (.)_{l}}{A_{l}} + \frac{(.)_{g}}{A_{g}},
\end{equation}
and furthermore 
\begin{equation}
\xi = \sqrt{ \rho^{*} (\rho_{l} - \rho_{g}) g_{n} \dd{h}{A_{l}} - \frac{\rho_{g} \rho_{l}}{A_{g} A_{l}} (u_g -u_l)^2}.
\end{equation}
The characteristic polynomial corresponding to the eigenvalue problem is given by:
\begin{equation}
\rho_{l} A_{g} (u_{l} - \lambda)^2 + \rho_{g} A_{l} (u_{g} - \lambda)^2 + (\rho_{l} - \rho_{g}) A_{g} A_{l} g_{n} \dd{h}{A_{l}} = 0,
\end{equation}
and the two eigenvalues $\lambda_{1,2}$ are real provided that:
\begin{equation}\label{eqn:IKH_limit}
(u_{g} - u_{l})^2 \leq \frac{\rho_{l} - \rho_{g}}{\rho^{*}} g_{n} \dd{h}{A_{l}},
\end{equation}
which is known as the inviscid Kelvin-Helmholtz limit. The (left) eigenvectors related to $\lambda_{1,2}$ are found as follows
\begin{align}
\vt{l}_{1} = \begin{pmatrix}
- \frac{\rho_{l}}{A_{l}} \frac{ \frac{\rho_{g}}{A_{g}} (u_{g} - u_{l}) - \xi}{\rho^{*}} & 
- \frac{\rho_{g}}{A_{g}} \frac{ \frac{\rho_{l}}{A_{l}} (u_{g} - u_{l}) + \xi}{\rho^{*}} &
-1 &
1
\end{pmatrix}, \label{eqn:eigenvector1}\\
\vt{l}_{2} = \begin{pmatrix}
- \frac{\rho_{l}}{A_{l}} \frac{ \frac{\rho_{g}}{A_{g}} (u_{g} - u_{l}) + \xi}{\rho^{*}} &
- \frac{\rho_{g}}{A_{g}} \frac{ \frac{\rho_{l}}{A_{l}} (u_{g} - u_{l}) - \xi}{\rho^{*}} &
-1 &
1
\end{pmatrix}.\label{eqn:eigenvector2}
\end{align}
The Riemann invariant related to each eigenvector follows by multiplying the original PDEs with the eigenvector. For $\lambda_{1}$ and $\vt{l}_{1}$ this yields
\begin{equation}\label{eqn:Riemann_invariant}
\vt{l}_{1} \cdot (\vt{A} \dd{\vt{W}}{t} + \vt{B} \dd{\vt{W}}{s} + \vt{S} ) = 0,
\end{equation}
giving the Riemann invariant
\begin{equation}
\xi \dd{A_{l}}{t} - \rho_{l} \dd{u_{l}}{t} + \rho_{g} \dd{u_{g}}{t} + \lambda_{1} \xi \dd{A_{l}}{s} - \lambda_{1} \rho_{l} \dd{u_{l}}{s} + \lambda_{1} \rho_{g} \dd{u_{g}}{s} + \frac{S_{l}}{A_{l}} -  \frac{S_{g}}{A_{g}} = 0,
\end{equation}
which can be written as
\begin{equation}
\xi \frac{\rd A_{l}}{\rd t}  - \rho_{l} \frac{\rd u_{l}}{\rd t} +  \rho_{g} \frac{\rd u_{g}}{\rd t} =  \frac{S_{g}}{A_{g}}  - \frac{S_{l}}{A_{l}},
\end{equation}
with $\frac{\rd }{\rd t} = \dd{}{t} + \lambda_{1} \dd{}{s}$. For the second eigenvalue and eigenvector $\lambda_{2}$, $\vt{l}_{2}$ we similarly get
\begin{equation}
- \xi \frac{\rd A_{l}}{\rd t}  - \rho_{l} \frac{\rd u_{l}}{\rd t} +  \rho_{g} \frac{\rd u_{g}}{\rd t} = \frac{S_{g}}{A_{g}}  - \frac{S_{l}}{A_{l}},
\end{equation}
with $\frac{\rd }{\rd t} = \dd{}{t} + \lambda_{2} \dd{}{s}$. These eigenvalues and Riemann invariants are similar to the ones typically obtained for the two-equation system mentioned in section \ref{sec:alternative_forms}, see \cite{Akselsen2016, Guo2002, Picchi2016a}, and will be used for our characteristic boundary condition treatment in section \ref{sec:bc}. These two eigenvalues and Riemann invariants reflect the hyperbolic character of the two-fluid model, but the effect of the third and fourth eigenvalue has been lost, and in particular the elliptic behaviour of the pressure is missing. 

Here we extend the characteristic analysis to include the third and fourth eigenvalue of the system. This requires the consideration of the \textit{inverse} eigenvalue problem, $\vt{A} \vt{v} = \mu \vt{B} \vt{v}$, with the determinant equation
\begin{equation}
\text{det} (\mu \vt{B} - \vt{A}) = 0,
\end{equation}
leading to the following characteristic polynomial
\begin{equation}
\mu^2 \left( \rho_{l} A_{g} (\mu u_{l} -1)^2 + \rho_{g} A_{l} (\mu u_{g} - 1)^2) + \mu^2 (\rho_{l}-\rho_{g}) A_{g} A_{l} g_{n} \dd{h}{A_{l}} \right) =0.
\end{equation}
This equation has solutions $\mu_{1,2} = 1/\lambda_{1,2}$, and $\mu_{3,4} = 0$, indicating that $\lambda_{3,4}$ are infinite (see also Drew and Passman \cite{Drew1998}). The left eigenvector of the inverse eigenvector problem is the same as of the original eigenvector problem, and for $\mu_{3}=0$ the condition $\vt{l}_{3} \cdot \vt{A} = 0$ leads to:
\begin{equation}
\vt{l}_{3} = \begin{pmatrix}
1 & 1 & 0 & 0
\end{pmatrix}.
\end{equation}
The corresponding `Riemann invariant' follows as:
\begin{equation}
\vt{l}_{3} \cdot (\vt{A} \dd{\vt{W}}{t} + \vt{B} \dd{\vt{W}}{s} + \vt{S} ) = 0,
\end{equation}
leading to
\begin{equation}
(u_{l} - u_{g}) \dd{A_{l}}{s} + A_{l} \dd{u_{l}}{s} + A_{g} \dd{u_{g}}{s} = 0,
\end{equation}
or equivalently
\begin{equation}\label{eqn:mixture_velocity_constraint3}
\text{\underline{C1}:} \qquad \dd{}{s} \left( u_{l} A_{l} + u_{g} A_{g} \right) = 0.
\end{equation}
Thus, this means that: \textit{the third `Riemann invariant' corresponds to the volumetric flow constraint \eqref{eqn:mixture_velocity_constraint1}.}

To find an eigenvector associated with the fourth eigenvalue $\mu_{4}=0$, we resort to a generalized eigenvector via the application of Jordan theory:
\begin{equation}
\vt{l}_{4} \cdot (\vt{A} - \mu_{4} \vt{B}) = \vt{l}_{3} \cdot \vt{B} = \begin{pmatrix}
u_{l} - u_{g} & A_{l} & A_{g} & 0
\end{pmatrix},
\end{equation}
which leads to ($\mu_{4}=0$)
\begin{equation}
\vt{l}_{4} = \begin{pmatrix}
u_{l} & u_{g} & A_{l} / \rho_{l} & A_{g} / \rho_{g}
\end{pmatrix}.
\end{equation}
Carrying out again the multiplication with the original PDE leads to
\begin{equation}\label{eqn:constraint3}
\dd{}{t} (u_{l} A_{l} + u_{g} A_{g}) + \dd{}{s} (u_{l}^2 A_{l} + u_{g}^2 A_{g}) = - \left(\frac{A_{l}}{\rho_{l}} + \frac{A_{g}}{\rho_{g}}\right) \dd{p}{s} + \frac{E_{g}}{\rho_{g}} + \frac{E_{l}}{\rho_{l}} + \frac{S_{g}}{\rho_{g}}  + \frac{S_{l}}{\rho_{l}}.  
\end{equation}
This is a combination of the momentum equations (after division by the respective densities) which involves the time derivative of the total volumetric flow $V(t)$ (see \eqref{eqn:mixture_velocity_constraint2}):
\begin{equation}\label{eqn:constraint3_rewritten}
\left(\frac{A_{l}}{\rho_{l}} + \frac{A_{g}}{\rho_{g}}\right) \dd{p}{s}  =  -\frac{\rd V(t)}{\rd t} - \dd{}{s} (u_{l}^2 A_{l} + u_{g}^2 A_{g}) + \frac{E_{g}}{\rho_{g}} + \frac{E_{l}}{\rho_{l}} + \frac{S_{g}}{\rho_{g}}  + \frac{S_{l}}{\rho_{l}}.  
\end{equation}
Further differentiation with respect to $s$ gives the following pressure constraint:
\begin{equation}\label{eqn:Poisson}
\text{\underline{C2}:} \qquad \dd{}{s} \left( \left( \frac{A_{g}}{\rho_{g}} + \frac{A_{l}}{\rho_{l}} \right) \frac{\partial p }{\partial s} \right) = - \frac{\partial^2 }{\partial s^2} (u_{g}^2 A_{g} + u_{l}^2 A_{l} ) + \dd{}{s} \left( \frac{E_{g}}{\rho_{g}} + \frac{E_{l}}{\rho_{l}} \right) + \dd{}{s}\left( \frac{S_{g}}{\rho_{g}} + \frac{S_{l}}{\rho_{l}} \right).
\end{equation}
\textit{The fourth `Riemann invariant' corresponds to a Poisson-type equation for the pressure}. Effectively, this equation is obtained by taking the time derivative of the mixture equation \eqref{eqn:mixture_velocity_constraint3}, and by substituting the spatially-differentiated momentum equations. It hinges on the fact that \textit{the time derivative of the spatial terms in the mass equations is equal to the spatial derivative of the temporal terms in the momentum equations}. This is an important observation that will be used in constructing a discrete approximation and a consistent discrete pressure equation. Equations 
\eqref{eqn:mixture_velocity_constraint3} and \eqref{eqn:Poisson} are two hidden constraints of the continuous model equations and can be used for constructing an efficient numerical solution algorithm.

Strictly speaking, the full two-fluid model is parabolic, because it has real eigenvalues with a degenerate set of eigenvectors. Practically speaking, the system of equations has a hyperbolic part, with two real eigenvalues, and an elliptic part, corresponding to the Poisson-like equation for the pressure. The pressure is a Lagrange multiplier that makes the mixture velocity field divergence free. 

To conclude, the first novel result of this paper is that \textit{the two `Riemann invariants' associated with the infinite eigenvalues are the hidden constraints of the system. They are satisfied instantaneously, at each moment in time.} This means, for example, that no explicit initial condition for the pressure is required, since it follows from the initial condition for the phase velocities and hold-up fractions. Furthermore, the presence of infinite eigenvalues requires an (at least partially) implicit time integration strategy. 

The insights obtained from the characteristic analysis in this section will be employed in the next sections to construct a new boundary condition treatment and to interpret our new time integration strategy. 

\section{Spatial discretization and boundary conditions}\label{sec:spatial_discretization_bc}

\subsection{Discretization at interior points}\label{sec:spatial_discretization}
The spatial discretization is on a staggered grid, consisting of $N$ `pressure' and $N+1$ `velocity' volumes. The midpoints of the velocity volumes lie on the faces of the pressure volumes. The pressure, density, hold-up and mass are defined in the centre of the pressure volumes, whereas the velocity and momentum are defined in the centre of the velocity volumes. For details we refer to \cite{Sanderse2017}. The unknowns are the vector of conservative variables $U(t)$:
\begin{equation}
U(t) = 
\begin{pmatrix}
m_{g} \\
m_{l} \\
I_{g} \\
I_{l}
\end{pmatrix} =
\begin{pmatrix}
[(\rho_{g} A_{g})_{1} \ldots (\rho_{g} A_{g})_{N}]^T \\
[(\rho_{l} A_{l})_{1} \ldots (\rho_{l} A_{l})_{N}]^T \\
[(\rho_{g} A_{g} u_{g})_{1/2} \ldots (\rho_{g} A_{g} u_{g})_{N+1/2}]^T \\
[(\rho_{l} A_{l} u_{l})_{1/2} \ldots (\rho_{l} A_{l} u_{l})_{N+1/2}]^T 
\end{pmatrix},
\end{equation}
and the pressure:
\begin{equation}
p(t) = [p_{1} \ldots p_{N}]^T.
\end{equation}
Note that $m_{g}$, $m_{l}$, $I_{g}$ and $I_{l}$ are vectors, containing mass and momentum (per unit pipe length) at the pressure and velocity volumes, respectively. $U(t)$ and $p(t)$ are both a function of time only. 

\begingroup
\begin{figure}[hbtp]
\fontsize{10pt}{12pt}\selectfont
\centering 
\def\svgwidth{ \textwidth}
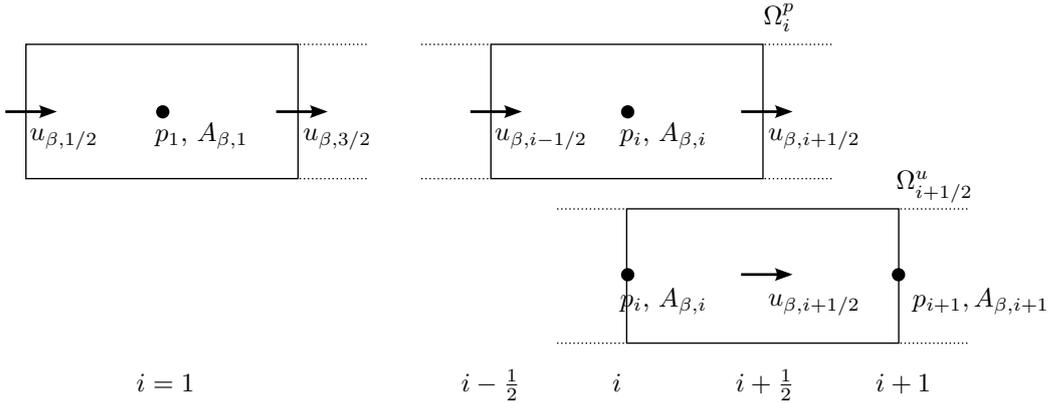 
\caption{Staggered grid layout including left boundary.}
\label{fig:GridLayout}
\end{figure}
\endgroup


We start with conservation of mass for phase $\beta$ ($\beta$ is liquid or gas).
Integration of equation \eqref{eqn:conservation_mass} in $s$-direction over a pressure volume gives:
\begin{equation}
\frac{\rd}{\rd t}\left(m_{\beta,i} \Delta s_{i} \right) + I_{\beta,i+1/2} -I_{\beta,i-1/2} = 0,\label{eqn:mass_semidiscrete}
\end{equation}
where the convective fluxes can be directly expressed in terms of the momenta $I$ at the staggered locations, so no approximation is involved in this term. 

For conservation of momentum we proceed in a similar way. Integration of \eqref{eqn:conservation_momentum} in $s$-direction over a velocity volume gives:
\begin{multline}
\frac{\rd}{\rd t}\left( I_{\beta,i+1/2} \Delta s_{i+1/2} \right) + \left(\rho_{\beta} A_{\beta}\right)_{i+1} (u_{\beta,i+1})^2 - \left(\rho_{\beta} A_{\beta} \right)_{i} (u_{\beta,i})^2 = -A_{\beta,i+1/2}\left(p_{i+1} - p_{i} \right) \\ + E_{\beta,i+1/2} + S_{\beta,i+1/2} \Delta s_{i+1/2}, \label{eqn:momentum_semidiscrete} 
\end{multline}
and the level gradient terms for the gas and liquid are given by (+ for gas, - for liquid)
\begin{equation}
E_{\beta,i+1/2} = \rho_{\beta} g_{n} \left( \left(h A_{\beta} \pm \frac{1}{12} P_{gl}^3 \right)_{i+1} - \left(h A_{\beta} \pm \frac{1}{12} P_{gl}^3 \right)_{i} \right). 
\end{equation}
The convective term in the momentum equation requires approximation; in the test cases in this work we have simply used a central approximation, $u_{\beta,i}=\frac{1}{2}(u_{\beta,i+1/2} + u_{\beta,i-1/2})$.

The system is closed with the volume constraint \eqref{eqn:volume_constraint}, which is written in terms of the phase masses $m_{\beta}$  as: 
\begin{equation}
\frac{m_{g}}{\rho_{g}} + \frac{m_{l}}{\rho_{l}} - A = 0.
\end{equation}

We stress that the unsteady term in the momentum equation and the mass fluxes in the mass equations both contain the same quantity $I_{\beta} = \rho_{\beta} A_{\beta} u_{\beta}$. This ensures a discrete coupling between the mass and momentum equations in the same way as in the continuous case, where the incompressible pressure equation was derived by equating the time differentiation of the flux terms in the mass equations to the spatial differentiation of the unsteady terms in the momentum equations. This coupling is naturally achieved by the use of a staggered grid.

\subsection{Boundary conditions}\label{sec:bc}
The boundary conditions should be consistent with the characteristic directions \eqref{eqn:roots} and eigenvectors \eqref{eqn:eigenvector1}-\eqref{eqn:eigenvector2}.
We rewrite the Riemann invariants, similar to \eqref{eqn:Riemann_invariant}, in matrix notation as
\begin{equation}
\vt{L} \vt{A} \dd{\vt{W}}{t} + \vt{L} \vt{B} \dd{\vt{W}}{t} + \vt{L} \vt{S} = 0,
\end{equation}
where $\vt{L}$ is the matrix that contains the left eigenvectors $\vt{l}_{i}$ as rows. Written in full, this gives
\begin{multline}\label{eqn:Riemann_matrix}
\begin{pmatrix}
\xi & -\rho_{l} & \rho_{g} & 0 \\
-\xi & -\rho_{l} & \rho_{g} & 0 \\
0 & 0 & 0 & 0 \\
u_{l} - u_{g} & A_{l} & A_{g} & 0 
\end{pmatrix}
\begin{pmatrix}
\dd{A_{l}}{t} \\
\dd{u_{l}}{t} \\
\dd{u_{g}}{t} \\
\dd{p}{t}
\end{pmatrix}
+ \\
\begin{pmatrix}
\lambda_{1} \xi & -\lambda_{1} \rho_{l} & \lambda_{1} \rho_{g} & 0 \\
-\lambda_{2} \xi & -\lambda_{2} \rho_{l} &\lambda_{2}  \rho_{g} & 0 \\
u_{l} - u_{g} & A_{l} & A_{g} & 0 \\
u_{l}^2 - u_{g}^2 + \tilde{K} & 2 u_{l} A_{l} & 2 u_{g} A_{g} & \frac{A_{g}}{\rho_{g}} +  \frac{A_{l}}{\rho_{l}}
\end{pmatrix}
\begin{pmatrix}
\dd{A_{l}}{s} \\
\dd{u_{l}}{s} \\
\dd{u_{g}}{s} \\
\dd{p}{s}
\end{pmatrix}
+
\begin{pmatrix}
\frac{S_{l}}{A_{l}} - \frac{S_{g}}{A_{g}} \\
\frac{S_{l}}{A_{l}} - \frac{S_{g}}{A_{g}} \\
0 \\
-\frac{S_{g}}{\rho_{g}} - \frac{S_{l}}{\rho_{l}}\\
\end{pmatrix}
=0,
\end{multline}
where $\tilde{K} = - \frac{1}{\rho_{l}} \dd{K_{l}}{A_{l}} + \frac{1}{\rho_{g}} \dd{K_{g}}{A_{g}}$.

In classic hyperbolic systems, boundary conditions are prescribed corresponding to the number of incoming waves and the equations for the evolution of the remaining components follows from the characteristic equations, see e.g.\ \cite{Thompson1987}. For the compressible two-fluid model such an approach has been outlined in Olsen \cite{Olsen2004}, using an explicit time integration method. For the incompressible two-fluid model, the mixed hyperbolic-elliptic character of the equations requires an (at least partially) implicit approach, which we propose next.

We make the following observations based on \eqref{eqn:Riemann_matrix}. First, the two characteristic (hyperbolic) equations link changes in the hold-up $A_{l}$ to changes in the phase velocities $u_{l}$, $u_{g}$, independent of the pressure. Second, there is no evolution equation for the pressure, but only an equation for its spatial derivative (fourth row), which can be interpreted as a boundary condition for the pressure. This suggests the following strategy: update the hold-up and phase velocities according to a characteristic treatment (rows 1+2), and then solve a pressure equation with an (implied) boundary condition for the pressure (row 4), such that the mixture velocity field is divergence free (row 3).

%

\subsubsection{Inflow conditions and solid wall conditions}
At a pipe inlet, it is common that the liquid and gas mass flows are given as a function of time:
\begin{equation}
I_{g} (s=0,t) = \rho_{g} A_{g} u_{g} = I_{g,\text{inlet}}(t), \qquad I_{l}(s=0,t) =  \rho_{l} A_{l} u_{l} = I_{l,\text{inlet}}(t).
\end{equation}
We will assume that not only the actual mass flows are available, but also the time derivatives $\dot{I}_{g,\text{inlet}}$ and $\dot{I}_{l,\text{inlet}}$. Solid walls are a special case for which we impose the conditions
\begin{equation}
u_{g} = u_{l} = 0, \qquad \dd{u_{g}}{t} = \dd{u_{l}}{t} = 0.
\end{equation}

The challenge is to find an evolution equation for the hold-up $A_{l}$ at the boundary that is consistent with the characteristic equations.
\begingroup
\begin{figure}[hbtp]
\fontfamily{lmss}
\fontsize{10pt}{12pt}\selectfont
\centering 
\def\svgwidth{0.7\textwidth}
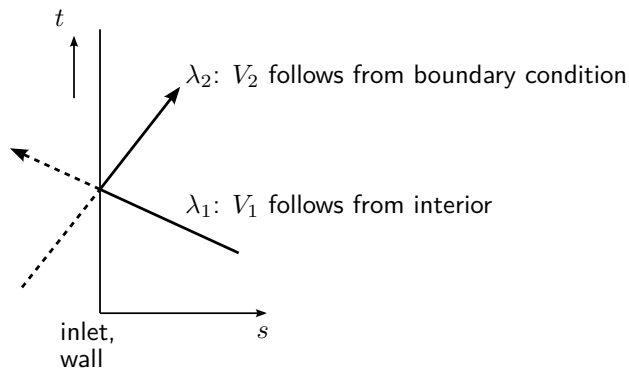 
\caption{Characteristics at the left boundary, used for the inflow and solid wall boundary condition.}
\label{fig:bcwall}
\end{figure}
\endgroup

Consider a boundary at $s=0$, as shown in figure \ref{fig:bcwall}, and assume an incoming wave $\lambda_{2}>0$ and an outgoing wave $\lambda_{1}<0$. The characteristic equations are
\begin{equation}
\begin{pmatrix}
\xi & -\rho_{l} & \rho_{g}  \\
-\xi & -\rho_{l} & \rho_{g}
\end{pmatrix}
\begin{pmatrix}
\dd{A_{l}}{t} \\
\dd{u_{l}}{t} \\
\dd{u_{g}}{t}
\end{pmatrix}
+
\begin{pmatrix}
\lambda_{1} \xi & -\lambda_{1} \rho_{l} & \lambda_{1} \rho_{g} \\
-\lambda_{2} \xi & -\lambda_{2} \rho_{l} &\lambda_{2}  \rho_{g} 
\end{pmatrix}
\begin{pmatrix}
\dd{A_{l}}{s} \\
\dd{u_{l}}{s} \\
\dd{u_{g}}{s}
\end{pmatrix}
+
\begin{pmatrix}
\frac{S_{l}}{A_{l}} - \frac{S_{g}}{A_{g}} \\
\frac{S_{l}}{A_{l}} - \frac{S_{g}}{A_{g}}
\end{pmatrix}
=0,
\end{equation}
or
\begin{equation}\label{eqn:characteristics}
\begin{pmatrix}
\xi & -\rho_{l} & \rho_{g} \\
-\xi & -\rho_{l} & \rho_{g}
\end{pmatrix}
\begin{pmatrix}
\dd{A_{l}}{t} \\
\dd{u_{l}}{t} \\
\dd{u_{g}}{t}
\end{pmatrix}
+
\begin{pmatrix}
\lambda_{1} V_{1}  \\
\lambda_{2} V_{2} 
\end{pmatrix}
+
\begin{pmatrix}
\frac{S_{l}}{A_{l}} - \frac{S_{g}}{A_{g}} \\
\frac{S_{l}}{A_{l}} - \frac{S_{g}}{A_{g}}
\end{pmatrix}
=0,
\end{equation}
where
\begin{equation}
\begin{pmatrix}
 V_{1}  \\
 V_{2} 
\end{pmatrix}
=
\begin{pmatrix}
 \xi & - \rho_{l} & \rho_{g} \\
- \xi & - \rho_{l} & \rho_{g}
\end{pmatrix}
\begin{pmatrix}
\dd{A_{l}}{s} \\
\dd{u_{l}}{s} \\
\dd{u_{g}}{s}
\end{pmatrix}.
\end{equation}

The characteristic equations feature time derivatives of the velocity; in order to obtain expressions in terms of mass flows we expand the mass flows in terms of velocity and hold-up changes as follows:
\begin{equation}
\dd{I_{g}}{t} = \rho_{g} A_{g} \dd{u_{g}}{t} + \rho_{g} u_{g} \dd{A_{g}}{t}.
\end{equation}
We can rewrite this equation into
\begin{equation}
\rho_{g} \dd{u_{g}}{t} = \frac{ \dd{I_{g}}{t} + \rho_{g} u_{g} \dd{A_{l}}{t} }{ A_{g}},
\end{equation}
and similarly
\begin{equation}
\rho_{l} \dd{u_{l}}{t} = \frac{ \dd{I_{l}}{t} - \rho_{l} u_{l} \dd{A_{l}}{t} }{ A_{l}}.
\end{equation}
Substituting these expressions into the characteristic equations \eqref{eqn:characteristics} yields
\begin{equation}
\begin{pmatrix}
\xi + \frac{\rho_{g} u_{g}}{A_{g}} + \frac{\rho_{l} u_{l}}{A_{l}} & \frac{-1}{A_{l}} & \frac{1}{A_{g}} \\
-\xi + \frac{\rho_{g} u_{g}}{A_{g}} + \frac{\rho_{l} u_{l}}{A_{l}} & \frac{-1}{A_{l}} & \frac{1}{A_{g}} \\
\end{pmatrix}
\begin{pmatrix}
\dd{A_{l}}{t} \\
\dd{I_{l}}{t} \\
\dd{I_{g}}{t}
\end{pmatrix}
+
\begin{pmatrix}
\lambda_{1} V_{1}  \\
\lambda_{2} V_{2} 
\end{pmatrix}
+
\begin{pmatrix}
\frac{S_{l}}{A_{l}} - \frac{S_{g}}{A_{g}} \\
\frac{S_{l}}{A_{l}} - \frac{S_{g}}{A_{g}}
\end{pmatrix}
=0.
\end{equation}

Since $\lambda_{1}$ corresponds to the wave that carries information from the interior to the boundary, $V_{1}$ is known from the solution in the interior of the domain. To obtain an equation for $V_{2}$ that uses the given boundary conditions for $\dd{I_{g}}{t}$ and $\dd{I_{l}}{t}$, we add the two characteristic equations in such a way that the hold-up term disappears:
\begin{equation}
2 \xi
\begin{pmatrix}
\frac{-1}{A_{l}} & \frac{1}{A_{g}}
\end{pmatrix}
\begin{pmatrix}
\dd{I_{l}}{t} \\
\dd{I_{g}}{t}
\end{pmatrix}
+
(\xi - k) \lambda_{1} V_{1}  + (\xi + k) \lambda_{2} V_{2} 
+ 2\xi \left(\frac{S_{l}}{A_{l}} - \frac{S_{g}}{A_{g}} \right) 
=0,
\end{equation}
where 
\begin{equation}
k = (\rho u)^* = \frac{\rho_{g} u_{g}}{A_{g}} + \frac{\rho_{l} u_{l}}{A_{l}},
\end{equation}
and $\lambda_{2} V_{2}$ follows as
\begin{equation}\label{eqn:l2V2}
 \lambda_{2} V_{2} = - \frac{1}{\xi +k} \left[ (\xi -k) \lambda_{1} V_{1}  +
2\xi  \begin{pmatrix}
\frac{-1}{A_{l}} & \frac{1}{A_{g}}
\end{pmatrix}
\begin{pmatrix}
\dd{I_{l}}{t} \\
\dd{I_{g}}{t}
\end{pmatrix}
+2 \xi \left(\frac{S_{l}}{A_{l}} - \frac{S_{g}}{A_{g}} \right) \right].
\end{equation}
The solid wall boundary condition is a special case of this equation obtained when setting $k=0$ and $\dd{I_{l}}{t} =\dd{I_{g}}{t}=0$. The boundary condition for $\dd{A_{l}}{t}$ follows from subtracting the two characteristic equations, which gives 
\begin{equation}\label{eqn:holdup_charbceqn}
\left( \dd{A_{l}}{t} \right)_{\text{inlet}} =  \frac{\lambda_{2} V_{2} - \lambda_{1} V_{1}}{2 \xi},
\end{equation}
and substituting equation \eqref{eqn:l2V2} for $\lambda_{2} V_{2}$. The treatment of an inflow or solid boundary at the right side of the domain ($s=L$) follows in a similar manner.

Given the time derivative of the liquid hold-up fraction, equation \eqref{eqn:holdup_charbceqn}, we now have a complete description for the evolution of the conservative variables $U$ at  the boundary:
\begin{align}
\frac{\rd m_{g,1/2}}{\rd t} &= - \rho_{g} \left( \dd{A_{l}}{t} \right)_{\text{inlet}}, \label{eqn:bc1}\\
\frac{\rd m_{l,1/2}}{\rd t} &= \rho_{l} \left( \dd{A_{l}}{t} \right)_{\text{inlet}}, \label{eqn:bc2}\\
I_{g,1/2} &= I_{g,\text{inlet}}, \label{eqn:bc3} \\
I_{l,1/2} &= I_{l,\text{inlet}}.\label{eqn:bc4} 
\end{align}
The prescription of $I_{g}$ and $I_{l}$ in this way will be denoted as boundary conditions in \textit{strong form}. The disadvantage of the strong form is that order reduction can appear in the time integration method, as will become clear in section \ref{sec:accuracy}. An alternative is to specify the boundary conditions in \textit{weak form} via the time derivatives, and to integrate these values in time:
\begin{align}
\frac{\rd I_{g,1/2}}{\rd t} &=\dot{I}_{g,\text{inlet}}, \label{eqn:bc5} \\
\frac{\rd I_{l,1/2}}{\rd t} &= \dot{I}_{l,\text{inlet}}.\label{eqn:bc6} 
\end{align}
However, this introduces a time integration error at the boundary points. Our preferred approach is the strong form, i.e.\ assume that $I_{g,\text{inlet}}(t)$ and $I_{l,\text{inlet}}(t)$ are known, obtain the time derivatives for evaluating \eqref{eqn:l2V2} by analytical or numerical differentiation, and then apply a custom-made Runge-Kutta method that prevents order reduction.

The pressure on the boundary, $p_{1/2}$, is obtained by linear extrapolation of the pressure at the first two interior points, and its value does not influence the solution of $U$ at the interior points. The pressure at the interior points, $p_{1},\ldots p_{N}$, is determined by the pressure Poisson equation. As will be detailed when discussing the time integration method, the pressure Poisson equation is such that the velocity field becomes divergence-free (C1), which in turn makes that the volume constraint is satisfied (C0). Similar to the case of the (single-phase, 2D or 3D) incompressible Navier-Stokes equations \cite{Gresho1998,Veldman1990}, also in this 1D incompressible two-fluid model \textit{no boundary condition for the pressure} needs to be prescribed. The pressure boundary conditions are \textit{implied} by the discretization of the constraint equation, and are a consistent approximation to $\dd{p}{s}$ as given by the last equation in \eqref{eqn:Riemann_matrix}. Only in the case of an outflow boundary, at which a pressure value is prescribed, a pressure boundary condition is necessary. This is detailed next.


\subsubsection{Outflow conditions}
At outflow boundaries a similar approach is taken as is common for the (single-phase) incompressible Navier-Stokes equations \cite{Sanderse2013}. We consider a `half' finite volume for the momentum equations, as shown in figure \ref{fig:halfvol}. In this volume we solve the momentum equations in the same way as in the interior points, with the only difference that the pressure is specified in the boundary point: $p_{N+1/2} = p_{\text{outlet}}$. The momentum equations yield expressions for the time derivatives of $I_{g,N+1/2}$ and $I_{l,N+1/2}$, which are subsequently used together with the characteristic treatment outlined before to obtain an equation for the hold-ups, viz. \eqref{eqn:holdup_charbceqn}. 

\begingroup
\begin{figure}[hbtp]
\fontsize{10pt}{12pt}\selectfont
\centering 
\def\svgwidth{0.7 \textwidth}
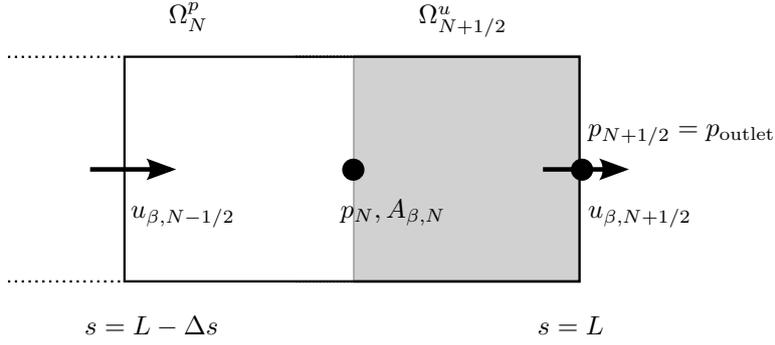 
\caption{Velocity volume (shaded) at an outflow boundary}
\label{fig:halfvol}
\end{figure}
\endgroup

\subsection{Summary}
In this section the characteristic directions and corresponding Riemann invariants derived in section \ref{sec:governing_equations} have been used to derive the second novel contribution of this paper: a new boundary condition treatment. The treatment is not only consistent with the wave directions but also with the constraints through the use of the pressure equation. No boundary conditions need to be prescribed for this pressure equation, which is consistent with the continuous model formulation for which no pressure boundary conditions were specified either.

\section{Constraint analysis of differential algebraic equations}\label{sec:temporal_properties}
In this section we treat the constraint analysis of the semi-discrete equations as indicated in the second row in figure \ref{fig:DAE_ODE}.
\subsection{DAE formulation}
The semi-discrete equations \eqref{eqn:mass_semidiscrete}-\eqref{eqn:momentum_semidiscrete} supplemented with boundary conditions \eqref{eqn:bc1}-\eqref{eqn:bc2} can be rewritten in terms of the differential unknowns $U(t) = [m(t),I(t)]^T$ and algebraic unknowns $p(t)$ as a semi-explicit differential-algebraic equation (DAE) system:
\begin{align}
\dot{U}(t) := \frac{\rd U(t)}{\rd t} &= F(U(t),p(t),t),\label{eqn:semi_discrete1}\\
\text{\underline{C0}:} \qquad g(U(t)) &= 0,\label{eqn:semi_discrete2}
\end{align}
where 
\begin{equation}
F(U(t),p(t)) := \begin{bmatrix} F_{m} (I(t),t) \\ F_{I} (m(t),I(t),t) - H (m(t)) p(t) \end{bmatrix}, \label{eqn:F}
\end{equation}
and
\begin{align}
F_m (I(t),t) &:= - D I(t) - b(t), \label{eqn:Fm} \\
 g(U(t)) &:= Q m (t) - A. \label{eqn:semidiscrete_g}
\end{align}
It is important to note that - in contrast to the compressible two-fluid model - the constraint equation does not depend on the pressure. This has important consequences for the index of the DAE system.

We have defined the following variables:
\begin{align}
&m(t) = \begin{bmatrix} m_{g}(t) \\ m_{l}(t) \end{bmatrix}, \qquad I(t) = \begin{bmatrix} I_{g}(t) \\ I_{l}(t) \end{bmatrix}, \qquad b(t)= \begin{bmatrix} b_{g}(t) \\ b_{l}(t) \end{bmatrix}, 
\end{align}
and the matrices
\begin{align}
& H(m(t)) = \begin{bmatrix} R m_g(t) \\ R m_l(t) \end{bmatrix} G, & Q = \begin{bmatrix} \mathbb{I}/\rho_{g} & \mathbb{I}/\rho_{l} \end{bmatrix}, \qquad D = \begin{bmatrix} D_p & 0 \\ 0 & D_p \end{bmatrix}. 
\end{align}
$\mathbb{I}$ is the $N\times N$ identity matrix, and $D_p$ and $G$ are differencing operators that compute differences from volume faces to midpoints or vice versa. For example, with periodic boundary conditions we have the $N\times N$ matrices 
\begin{equation}
D_p = \begin{pmatrix}
-1 & 1 & &  \\
  & \ddots & \ddots &  \\
  & &   -1 & 1 \\
1 &  &  & -1
\end{pmatrix}, \qquad
G = \begin{pmatrix}
1 &  &  & -1 \\
-1 & 1 & &  \\
  & \ddots & \ddots &  \\
  & &   -1 & 1 \\
\end{pmatrix}.
\end{equation}
$R$ is an interpolation matrix that computes averages on cell interfaces based on cell midpoint values. $b_{\beta}$ contains boundary conditions, e.g. $I_{g,\text{inlet}}$, possibly depending on time but not on the solution.

In what follows, the dependence of $F$ and $g$ on $U$ and $p$, and their dependence on $(t)$ will be omitted in the equations. Note that explicit time-dependence, e.g.\ due to time-dependent source terms or time-dependent boundary conditions, can be accommodated by adding time as an  unknown to $U$ and adding the equation $\dot{t}=1$.

\subsection{Derivation of the index and hidden constraints}
The index of the DAE system is an important concept that gives the theoretical framework for studying the order conditions and the associated order of accuracy of the time integration methods that will be proposed in section \ref{sec:time_integration}. The index and any hidden constraints are revealed by differentiating the DAE system in time:
\begin{align}
\ddot{U} &= \frac{\partial F}{\partial U} \dot{U} + \frac{\partial F}{\partial p} \dot{p}, \label{eqn:DAE_incomp1_diff1} \\
\frac{\partial g}{\partial U} \dot{U} &= 0. \label{eqn:DAE_incomp2_diff1}
\end{align}
The differentiated constraint can be rewritten after substituting equation \eqref{eqn:semi_discrete1}: 
\begin{equation}\label{eqn:PPE}
\frac{\partial g}{\partial U} F(U,p) = 0.
\end{equation}
Evaluating this equation by using equations \eqref{eqn:semi_discrete2} and \eqref{eqn:F} gives
\begin{equation}
\begin{split}
\frac{\partial g}{\partial U} F(U,p) &  =
\begin{pmatrix} 
\frac{\partial g}{\partial m} & \frac{\partial g}{\partial I}  
\end{pmatrix}
\begin{pmatrix} 
F_{m} \\ F_{I}
\end{pmatrix}  
= Q F_{m} = 0, 
\end{split}
\end{equation}
and after substituting equation \eqref{eqn:Fm} we obtain:
\begin{equation}\label{eqn:hidden_constraint1}
\text{\underline{C1}}: \qquad \boxed{M I + r = 0,}
\end{equation}
with $M = QD$ and $r=Q b$. For example, with inflow conditions at $s=0$, and outflow conditions at $s=L$, we have
\begin{equation}
r (t) = \begin{pmatrix}
- I_{g,\text{inlet}}(t)/ \rho_{g} -  I_{l,\text{inlet}}(t)/ \rho_{l} \\
0 \\
\vdots \\
0
\end{pmatrix} = \begin{pmatrix}
-V(t) \\
0 \\
\vdots \\
0
\end{pmatrix}.
\end{equation}
Equation \eqref{eqn:hidden_constraint1} is the semi-discrete equivalent of the volumetric flow constraint that was encountered as a Riemann invariant in the characteristic analysis, viz.\ equation  \eqref{eqn:mixture_velocity_constraint3}. 

In contrast to the compressible two-fluid model, a single differentiation does not yield an equation for the pressure, because $F_m$ is not a function of $p$, i.e.\ $\frac{\partial g}{\partial U} \frac{\partial F}{\partial p} = 0$. 
To obtain a pressure equation we need to further differentiate the constraint:
\begin{align}
\frac{\partial g}{\partial U} \ddot{U} &= 0. \label{eqn:DAE_incomp2_diff2}
\end{align}
Note that $\frac{\partial^2 g}{\partial U^2}=0$, since $Q$ and $A$ in equation \eqref{eqn:semidiscrete_g} are independent of $U$ and $t$. The constraint \eqref{eqn:DAE_incomp2_diff2} is rewritten by substituting \eqref{eqn:DAE_incomp1_diff1}, leading to
\begin{equation}
\frac{\partial g}{\partial U} \frac{\partial F}{\partial U} F + \frac{\partial g}{\partial U} \frac{\partial F}{\partial p}  \dot{p} = \frac{\partial g}{\partial U} \frac{\partial F}{\partial U} F = 0.\label{eqn:constraint2}
\end{equation}
This equation is the second hidden constraint (C2), the discrete equivalent of the constraint given by equation \eqref{eqn:constraint3}.
This constraint can be written in a more convenient form by substituting the expression for $F_{m}$, yielding  
\begin{equation}
\frac{\partial g}{\partial U} \frac{\partial F}{\partial U} F = 
- Q (D (F_{I}- H p)  + \dot{b}) = - M F_{I} + M H p  - \dot{r} =0,
\end{equation}
or
\begin{equation}\label{eqn:poisson}
\text{\underline{C2}}: \qquad \boxed{L p = M F_{I} + \dot{r}.}
\end{equation}
where $L = M H$ is a Laplace operator, being the discrete version of $\dd{}{s} ( (\frac{A_{g}}{\rho_{g}} + \frac{A_{l}}{\rho_{l}}) \dd{p}{s})$. $L$ is symmetric negative definite since the volume fractions are always non-negative. Equation \eqref{eqn:poisson} is the Poisson equation for the pressure in incompressible 1D multiphase flow and it corresponds to the pressure constraint C2 derived earlier in equation \eqref{eqn:Poisson}. 
A further differentiation is required to derive an ODE for the pressure, and to obtain the index of the DAE system:
\begin{align}
\frac{\partial g}{\partial U} \dddot{U} = 0.
\end{align}
Substituting the equation for $\dddot{U}$ (not shown here) into this equation gives an ODE for the pressure, because $\dd{g}{U} \dd{F}{U} \dd{F}{p}$ is non-singular: it is the Laplace operator $L$. Equivalently, one can take the time derivative of equation \eqref{eqn:poisson}. Since three differentiations were required to arrive at the ODE for the pressure, the index of the DAE system \eqref{eqn:semi_discrete1}-\eqref{eqn:semi_discrete2} is 3. To be precise, the DAE system is in semi-explicit form and has Hessenberg index 3 \footnote{Hessenberg form: ODEs coupled with constraint, with explicit identification of algebraic and differential variables, where the algebraic variables may all be eliminated using the same number of differentiations \cite{Ascher1998}.}. The existence of solutions for our index 3 problem is guaranteed since $F_{I}$ is \textit{linear} in $p$ \cite{Ostermann1993}.



\subsection{Summary}
The previous analysis highlights the third novel result of this paper: the semi-discrete pressure equation for incompressible flow can be derived in a structured manner by considering the framework of differential-algebraic equations. The system has index 3 and has the same constraints, C1 and C2 (for the volumetric flow and for the pressure, respectively), that were derived for the continuous equations in section \ref{sec:governing_equations}. 

As graphically shown in the second row of figure \ref{fig:DAE_ODE}, the index 3 system \eqref{eqn:semi_discrete1}-\eqref{eqn:semi_discrete2} can be rewritten in terms of an index 2 system (equations \eqref{eqn:semi_discrete1} and \eqref{eqn:hidden_constraint1}), or an index 1 system (equations \eqref{eqn:semi_discrete1} and \eqref{eqn:poisson}). This index reduction can be useful in order to apply existing (index 1 or 2) time integration methods. However, care has to be taken that the original constraint C0, \eqref{eqn:semidiscrete_g}, is satisfied. In the next section we propose time integration methods based on the index 3 formulation that preserve all three constraints C0, C1 and C2.


\section{New time integration methods}\label{sec:time_integration}

\subsection{Explicit Runge-Kutta methods}
In this section we focus on the third row of figure \ref{fig:DAE_ODE}: discretization of equations \eqref{eqn:semi_discrete1}-\eqref{eqn:semi_discrete2} in time with an explicit Runge-Kutta method. Application of an explicit RK-method is not trivial due to the presence of the constraint, which has to be satisfied at the new stage level, and consequently introduces a degree of implicitness to the system. Related to this issue is the evaluation of the pressure term in the momentum equation. A possible approach is to discretize the pressure equation corresponding to the pressure constraint C2, \eqref{eqn:poisson}, but this does not directly guarantee compatibility with the Runge-Kutta time integration strategy and the constraint. Instead, we will take the following approach: first, discretize the system of equations with a so-callled half-explicit Runge-Kutta method \cite{Hairer1989}, and then perform substitutions on the fully discrete system, such that a fully discrete pressure equation consistent with constraints C0 and C1 will be obtained. 

First, we assume that the initial conditions are consistent and that they fulfil the hidden constraints, i.e. 
\begin{align}
& \text{\underline{C0}:} & Q m^{0} &= A, \label{eqn:IC1}\\
& \text{\underline{C1}:} & M I^{0} + r(t^0) &= 0, \label{eqn:IC2}\\
& \text{\underline{C2}:} & L^{0} p^{0} &= M F_{I}^{0} + \dot{r}(t^0). \label{eqn:IC3}
\end{align}
The half-explicit Runge-Kutta method for the index 3 system then follows by advancing the differential part of the equations with an explicit procedure, but requiring the constraint to be satisfied at each stage of the Runge-Kutta method:
\begin{align}
m_{i} & = m^{n} + \Delta t \sum_{j=1}^{i-1} a_{ij} F_{m,j}, \label{eqn:RK_mass}\\
I_{i} & = I^{n} + \Delta t \sum_{j=1}^{i-1} a_{ij} (F_{I,j} - H_{j} p_{j}), \label{eqn:RK_momentum}\\
Q m_{i} &= A,\label{eqn:RK_constraint}
\end{align}
followed by an update to the new time level $t^{n+1}$ by combining the stage values, again constrained by the volume constraint:
\begin{align}
m^{n+1} & = m^{n} + \Delta t \sum_{i=1}^{s} b_{i} F_{m,i}, \label{eqn:RK_mass_new} \\
I^{n+1} & = I^{n} + \Delta t \sum_{i=1}^{s} b_{i} (F_{I,i} - H_{i} p_{i}), \label{eqn:RK_momentum_new}\\
Q m^{n+1} &= A.\label{eqn:RK_constraint_new}
\end{align}
The subscript $i$ denotes the stage level, and should not be confused with the spatial index used earlier. $m_{i}$, $I_{i}$ and $p_{i}$ are approximations to $m(t_{i})$, $I(t_{i})$, $p(t_{i})$, where $t_{i} = t^n + c_{i} \Delta t$. The coefficients $a_{ij}$, $b_{i}$ and $c_{i}$ form the Butcher tableau that fully defines the Runge-Kutta method. 

\textit{We stress that equations \eqref{eqn:RK_mass}-\eqref{eqn:RK_constraint} and \eqref{eqn:RK_mass_new}-\eqref{eqn:RK_constraint_new} fully define the time integration method and hence determine its solution. The subsequent `substitution' steps presented next are simply a reformulation of these equations in order to arrive at a predictor-corrector type algorithm that involves the solution of a Poisson equation.}



The substitution process follows the same route as the semi-discrete case: substitute the mass equations \eqref{eqn:RK_mass} into the constraint \eqref{eqn:RK_constraint} to get:
\begin{equation}\label{eqn:RK_derivation1}
\Delta t \sum_{j=1}^{i-1} a_{ij} Q F_{m,j} = A - Q m^{n} = 0,
\end{equation}
which becomes, after substituting the expression $F_{m}=-D I - b$,
\begin{equation}\label{eqn:RK_derivation2}
 \sum_{j=1}^{i-1} a_{ij}  \left( M I_{j} + r_j \right)= 0.
\end{equation}
Since we use an explicit Runge-Kutta method, this reduces to
\begin{equation}\label{eqn:RK_derivation3}
\text{\underline{C1}:} \qquad M I_{i-1} + r_{i-1} = 0.
\end{equation}
In other words, \textit{the volume constraint at stage $i$ can be rewritten in terms of the volumetric flow constraint evaluated at stage $i-1$}. Equation \eqref{eqn:RK_derivation3} is the fully discrete equivalent of equations \eqref{eqn:mixture_velocity_constraint1} and \eqref{eqn:hidden_constraint1}.

Subsequently, the expression for $I_{i-1}$ follows from substituting the momentum equation, equation \eqref{eqn:RK_momentum}, and we obtain:
\begin{equation}
 M \left( I^{n} + \Delta t \sum_{j=1}^{i-2} a_{i-1,j} (F_{I,j} - H_{j} p_{j}) \right) + r_{i-1} = 0.
\end{equation}
Rewriting gives
\begin{equation}\label{eqn:RK_Poisson_lag}
\text{\underline{C2}:} \quad  a_{i-1,i-2} L_{i-2} p_{i-2} =  a_{i-1,i-2} M F_{I,i-2} +  \sum_{j=1}^{i-3} a_{i-1,j}  \left( M F_{I,j} - L_{j} p_{j})\right) + \frac{r_{i-1} - r^{n}}{\Delta t}.
\end{equation}
This is the Poisson equation for the pressure $p_{i-2}$, i.e.\ the second hidden constraint evaluated at stage $i-2$. In other words, \textit{the volume constraint at stage $i$ leads to the pressure  constraint evaluated at stage $i-2$}:
\begin{equation}
Q m_{i} = A \quad \rightarrow \quad  M I_{i-1} + r_{i-1} = 0 \quad \rightarrow \quad  L_{i-2} p_{i-2} = \ldots
\end{equation}
For example, in a 3-stage method, $p_{1}$ is the pressure that ensures that the velocity field $I_{2}$ is divergence free, which in turn makes that the hold-ups $m_{3}$ satisfy the volume constraint. This is indicated in figure \ref{fig:RKconstraints}. The index 3 nature of the system manifests itself in that \textit{two} substitution steps are required to obtain an equation for the pressure, in contrast to the single-phase, index 2 case, where only one substitution is required. This is consistent with the derivation of the pressure equation \eqref{eqn:poisson} and the observations in \cite{Ostermann1993}.


\begingroup
\begin{figure}[hbtp]
\centering 
\def\svgwidth{0.9\textwidth}
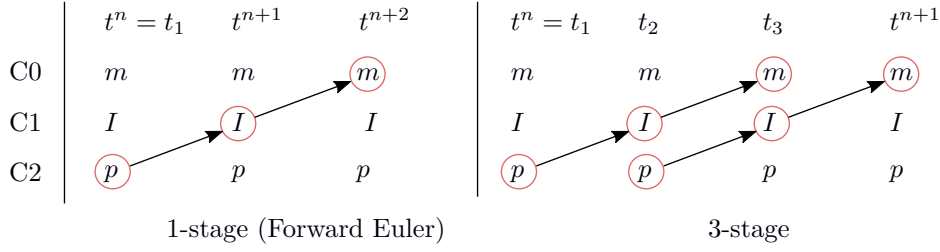 
\caption{Example of constraint progression for 1-stage and 3-stage Runge-Kutta method.\label{fig:RKconstraints}}
\end{figure}
\endgroup

The pressure equation \eqref{eqn:RK_Poisson_lag} can be rewritten in a more convenient form by simply changing the indices:
\begin{equation}\label{eqn:RK_poisson}
a_{i,i-1} L_{i-1} p_{i-1} =  a_{i,i-1} M F_{I,i-1} +  \sum_{j=1}^{i-2} a_{i,j}  \left( M F_{I,j} - L_{j} p_{j})\right) + \frac{r_{i} - r^{n}}{\Delta t}.
\end{equation}

With equation \eqref{eqn:RK_poisson} we rewrite \eqref{eqn:RK_mass}-\eqref{eqn:RK_constraint} into a `predictor-corrector' algorithm, similar to the fractional-step methods commonly used in incompressible Navier-Stokes algorithms. The first stage is trivial and gives $F_{I,1}=F_{I}(m_n,I_n,p_n)$ and $F_{m,1}=F_{m}(I_n)$. 
The subsequent stages are given by ($i=2, \ldots, s$)
\begin{empheq}[box=\fbox]{align}
m_{i} & = m^{n} + \Delta t \sum_{j=1}^{i-1} a_{ij} F_{m,j}, \\
I_{i}^{*} & = I^{n} + \Delta t \left( \sum_{j=1}^{i-1} a_{ij} F_{I,j} - \sum_{j=1}^{i-2} a_{ij} H_{j} p_{j} \right), \\
L_{i-1} \phi_{i-1} &= M I_{i}^{*} + r_{i}, \label{eqn:Poisson_phi}\\
I_{i} & = I_{i}^{*} - H_{i-1} \phi_{i-1}.
\end{empheq}
Here we defined $\phi_{i-1}$ such that $p_{i-1} = \phi_{i-1} / (a_{i,i-1} \Delta t)$. The final update to the new time level $n+1$ is given by
\begin{empheq}[box=\fbox]{align}
m^{n+1} & = m^{n} + \Delta t \sum_{i=1}^{s} b_{i} F_{m,i},\\
I^{n+1,*} & = I^{n} + \Delta t \left( \sum_{i=1}^{s} b_{i} F_{I,i} - \sum_{i=1}^{s-1} b_{i} H_{i} p_{i} \right), \\
L_{s} \phi_{s} &= M I^{n+1,*} + r^{n+1},\label{eqn:RK_pressure_last_stage}\\
I^{n+1} & = I^{n+1,*} - H_{s} \phi_{s}.
\end{empheq}

We note the following:
\begin{itemize}
\item The solution at the intermediate stages ($m_{i}$ and $I_{i}$) does not depend on the pressure at the start of the time step, $p^{n}$.
\item The predictor-corrector formulation does not contain a splitting error (see also \cite{Sanderse2012a}).
\item The Laplace operator $L_{i}=M H_{i}$ changes from stage to stage because the gradient operator $H_{i}$ depends on the solution.  An $s$-stage Runge-Kutta method requires $s$ Poisson solves per time step. 
\item It is necessary that the subdiagonal of the Butcher tableau, $a_{i+1,i}$, has all entries nonzero.
\item There is no start-up problem. The first stage is trivial; the second stage gives $m_{2} = m^{n} + a_{21} \Delta t F_{m} (I_{n})$ which automatically satisfies the volume constraint $Q m_{2} = Q m_{n} = A$, as long as the initial conditions are consistent. Subsequently, $\phi_{1}$ is determined such that $I_{2}$ satisfies the volumetric flow constraint $M I_{2} + r_{2} = 0$, and such that $Q m_{3} = Q m^{n} = A$.
\end{itemize}

\subsection{Accuracy of the Runge-Kutta method}\label{sec:accuracy}
The analysis of the temporal accuracy of Runge-Kutta methods applied to general index 3 DAE systems is not simple - see for example \cite{Jay1995,Ostermann1993}. Fortunately, because the volume constraint of the two-fluid model is simple, the fully discrete equations can be rewritten as if they were derived from an index 2 formulation. This is indicated in figure \ref{fig:DAE_ODE} by the arrow `accuracy analysis'. The corresponding index 2 DAE formulation is derived as follows.
First write the half-explicit Runge-Kutta method of equations \eqref{eqn:RK_mass}, \eqref{eqn:RK_momentum} and \eqref{eqn:RK_derivation3} as
\begin{align}
U_{i} &= U^{n} + \Delta t \sum_{j=1}^{i-1} a_{ij} \left( \hat{F}_{j} - \hat{H}_{j} p_{j} \right), \label{eqn:RK_index2form1}\\
\hat{M} U_{i} + r_{i} &= 0, \label{eqn:RK_index2form2}
\end{align}
where the extended operators $\hat{(.)}$ are defined as
\begin{equation}
\hat{F} (U) = \begin{bmatrix}
F_{m} (I) \\ F_{I} (m,I)
\end{bmatrix},
\qquad
\hat{M} = \begin{bmatrix}
0 & M
\end{bmatrix},
\qquad 
\hat{H} (U) = \begin{bmatrix}
0 \\ H (m)
\end{bmatrix}.
\end{equation}
The gradient operator that acts on the pressure can be written as
\begin{equation}
H(m) p = (R m) \odot (\hat{G} p) = \text{diag} (Rm) \hat{G} p,
\end{equation}
where
\begin{equation}
R m = \begin{pmatrix}
R_{g} & 0 \\
0 & R_{l}
\end{pmatrix}
\begin{pmatrix}
m_{g} \\ m_{l} 
\end{pmatrix},
\qquad 
\hat{G} = \begin{pmatrix}
G \\ G
\end{pmatrix}.
\end{equation}
The corresponding semi-discrete system is then recognized as
\begin{align}
\frac{\rd U}{\rd t} &= F(U,p) \label{eqn:index2_1} \\
\hat{g}(U) := \hat{M} U + r &= 0. \label{eqn:index2_2}
\end{align}
A crucial observation is that this semi-discrete system (and the corresponding half-explicit Runge-Kutta method) has the same form as the single-phase incompressible Navier-Stokes equations, for which the order conditions were shown in \cite{Sanderse2012a}. However, there is one important exception in the current multi-phase flow problem: \textit{the pressure gradient term depends on the solution} due to the presence of the hold-up fractions: $H=H(m)$. This dependence leads to additional order conditions for the differential component $U$ for methods of order 3 or higher. These conditions are shown in Table 3 in \cite{Sanderse2012a}. 

For order 3, there is one additional order condition, which reads 
\begin{equation}\label{eqn:additional_condition_3}
\sum_{i,j} b_{i} c_{i} \omega_{ij} c_{j+1}^2 = \frac{2}{3},
\end{equation}
where $\omega$ is the inverse of the shifted Butcher tableau $\tilde{a}$:
\begin{equation}
\tilde{a} = \begin{pmatrix}
a_{21} & 0\\
\vdots &\ddots \\
a_{s1} & & a_{s,s-1} & 0\\ 
b_{1} & \ldots & b_{s-1} & b_{s}
\end{pmatrix}.
\end{equation}
The differential associated to order condition \eqref{eqn:additional_condition_3} is given by condition 10 in \cite{Sanderse2012a}, which reads
\begin{equation}\label{eqn:additional_condition_3_differential}
F_{pU} (F, (-\hat{g}_{U} F_{p})^{-1} \hat{g}_{UU} (F,F)) = - (R F_{m} ) \odot \hat{G} L^{-1} \ddot{r},
\end{equation}
where the right hand side is obtained by inserting system \eqref{eqn:index2_1}-\eqref{eqn:index2_2}. The dependence of the pressure gradient term on the hold-up fractions manifests itself in the fact that $F_{pU}\neq0$, which is in contrast to single-phase problems, where generally $F_{pU}=0$.

Requiring order condition \eqref{eqn:additional_condition_3} to be satisfied (in addition to the four classic order conditions for a third order Runge-Kutta method \cite{Hairer2006}), leads to a one-parameter family of methods, with the following Butcher tableau ($c_{2}\neq0$, $c_{2} \neq \frac{2}{3}$, $c_{2} \neq 1$):

\begin{equation}\label{eqn:RK3_tableaux}
\def\arraystretch{1.25}
\begin{array}{c|ccc}
\multicolumn{4}{c}{\text{}}\\
0& 0 \\
c_{2} & c_{2}  \\
1 & 1 + \frac{1-c_{2}}{c_{2}(3 c_{2} - 2)} & - \frac{1-c_{2}}{c_{2}(3 c_{2} - 2)}  \\
\hline 
& \frac{1}{2} - \frac{1}{6 c_{2}} & \frac{1}{6 c_{2}(1-c_{2})} & \frac{2-3c_{2}}{6(1-c_{2})} \\
\multicolumn{4}{c}{\text{}}\\
\multicolumn{4}{c}{\text{RK3}}\\
\end{array}
\end{equation}
We propose to use the value of $c_{2} = \frac{1}{2}$, which is such that the $b$-coefficients are all positive and such that most order conditions associated to fourth order are also satisfied (leading to a small truncation error).

For order 4, in total six additional order conditions appear \cite{Sanderse2012a}. It is shown in \cite{Brasey1993} that these cannot be satisfied when employing a four-stage method that also should satisfy the classical order conditions. This means that \textit{an explicit four-stage, fourth-order Runge-Kutta method does not exist} for the differential-algebraic equations arising in one-dimensional multiphase flow problems. To achieve fourth order, a five-stage method is needed. An example of a five-stage, fourth-order method that satisfies all additional order conditions is the HEM4 method described in \cite{Brasey1993}. 

A crucial remark is in place here. Equation \eqref{eqn:additional_condition_3_differential} indicates that in case $\ddot{r}=0$, the additional order condition for third order disappears. It turns out that this is also true for fourth order (see Table 3 in \cite{Sanderse2012a}). The case $\ddot{r}=0$ appears when either 
\begin{itemize}
\item the boundary conditions are prescribed in weak form; or
\item the boundary conditions are prescribed in strong form but they do not depend on time. 
\end{itemize}
In both cases \textit{no additional order conditions are present, and classic Runge-Kutta methods can be used.}

The issue of order reduction and additional order conditions is even more prominent for the accuracy of the pressure. In fact, the pressure $p_{s}$ obtained from $\phi_{s}/(b_{s} \Delta t)$ is generally only a first-order accurate approximation to the pressure $p(t^{n+1})$. However, there is a simple way to make the accuracy of the pressure the same as of the differential variables, namely by solving equation \eqref{eqn:poisson}, given the solution for $m^{n+1}$ and $I^{n+1}$:
\begin{equation}\label{eqn:poisson_additional}
L^{n+1} p^{n+1} = M F_{I}^{n+1} + \dot{r}^{n+1}.
\end{equation}
This can be performed as a postprocessing step whenever an accurate pressure is required, since $p^{n+1}$ does not influence the solution in the next time step.

To conclude, the following three Runge-Kutta methods will be considered in our test cases:
\begin{equation}\label{eqn:RK_tableaux}
\def\arraystretch{1.25}
\begin{array}{c|c}
\multicolumn{2}{c}{\text{}}\\
\multicolumn{2}{c}{\text{}}\\
\multicolumn{2}{c}{\text{}}\\
c_{i} & a_{ij} \\
\hline 
 & b_{j}  \\
\multicolumn{2}{c}{\text{}}\\
\multicolumn{2}{c}{\text{Tableau}}\\
\end{array}
\qquad
\begin{array}{c|cc}
\multicolumn{3}{c}{\text{}}\\
\multicolumn{3}{c}{\text{}}\\
0& 0 \\
1 & 1  \\
\hline 
& \frac{1}{2} & \frac{1}{2} \\
\multicolumn{3}{c}{\text{}}\\
\multicolumn{3}{c}{\text{RK2}}\\
\end{array}
\qquad  
\begin{array}{c|ccc}
\multicolumn{4}{c}{\text{}}\\
0& 0 \\
\frac{1}{2} & \frac{1}{2}  \\
1 & -1 & 2 \\
\hline 
& \frac{1}{6} & \frac{2}{3} & \frac{1}{6} \\
\multicolumn{4}{c}{\text{}}\\
\multicolumn{4}{c}{\underline{\text{RK3 - proposed}}}\\
\end{array}
\qquad
\begin{array}{c|cccc}
0 & 0 \\
\frac{1}{2} & \frac{1}{2}\\
\frac{1}{2} & 0 & \frac{1}{2}  \\
1 & 0 & 0 & 1 \\
\hline 
 & \frac{1}{6} & \frac{1}{3} & \frac{1}{3} & \frac{1}{6} \\
 \multicolumn{5}{c}{\text{}}\\
\multicolumn{5}{c}{\text{RK4}}\\
\end{array}
\end{equation}
The second order method is the explicit midpoint method. The third order method is our proposed method: tableau \eqref{eqn:RK3_tableaux} with $c_{2}=\frac{1}{2}$. The fourth order method is the classic fourth order method, which for our DAE problem is fourth order accurate only provided that $\ddot{r}=0$; otherwise it is third order accurate. When necessary, we will compare these methods to the RK3-SSP (strong-stability preserving) method \cite{Gottlieb2001} and the five-stage fourth order method HEM4 from \cite{Brasey1993}.

The stability domains of these explicit Runge-Kutta methods can be found in many time integration textbooks, e.g.\ \cite{Butcher2003} (note that the treatment of the constraint via the pressure equation is fully implicit and does not affect the stability). For convection-dominated problems (for example in case of the inviscid model), we have shown in previous work that the eigenvalues of the semi-discrete equations lie on the imaginary axis \cite{Sanderse2017}. From a stability point of view, RK3 and RK4 are therefore to be preferred, because the stability domain of these methods contains a part of the imaginary axis.

\subsection{Eliminating constraint drift}
In the solution of the pressure Poisson equation, \eqref{eqn:Poisson_phi}, numerical errors are generally introduced, for example due to the tolerance setting of an iterative method. This can lead to errors in the constraints C0 and C1, equations \eqref{eqn:RK_derivation1} and \eqref{eqn:RK_derivation2}, which could potentially accumulate over time. We present an approach to prevent this, based on the ideas outlined in \cite{Hirt1967a}. The most important observation is that, when reformulating the constraint equation to obtain Poisson equation \eqref{eqn:RK_poisson}, one should \textit{not} substitute $Q m^{n} = A$ or $M I_{j} + r_{j}=0$, but instead leave these terms inside the equations.

When keeping these terms in the equations, a pressure Poisson equation that is similar to equation \eqref{eqn:Poisson_phi} follows, but with an additional term that involves the accuracy with which the constraints C0 and C1 have been satisfied. For the intermediate stages ($i = 2, \ldots, s$) the Poisson equation is
\begin{align}
L_{i-1} \phi_{i-1} &=  M I^{*}_{i} + r_{i} + \eta_{i}, \label{eqn:ppe_correction_stage} \\
\eta_{i} & = \left( \sum_{j=1}^{i-1} a_{i+1,j} (M I_{j} + r_{j}) - \frac{Q m^{n} - A}{\Delta t} \right) / a_{i+1,i}, 
\end{align}
and the Poisson equation for the update to the next time step is
\begin{align}
L_{s} \phi_{s} &=  M I^{n+1,*} + r^{n+1} + \eta^{n+1},\label{eqn:ppe_correction} \\
\eta^{n+1} &= - \frac{Q m^{n+1} - A}{\Delta t} / a_{21}.
\end{align}


As mentioned, constraint errors can be caused by the accuracy with which the Poisson equation is solved. In principle, the solution of the Poisson equations \eqref{eqn:ppe_correction_stage} and \eqref{eqn:ppe_correction} is straightforward, since the matrix $L$ is tri-diagonal and symmetric negative definite. This type of equation can be solved efficiently with a direct solver or with a preconditioned conjugate gradient solver. In case pressure boundary conditions are not involved, the pressure is determined up to a constant and consequently the matrix $L$ is singular. In principle, this is not a problem, as long as a constant solution lies in the null space of $L$. However, a direct solver can have difficulties with such a system and therefore we prefer the conjugate gradient solver, which has been used in the test cases reported here. In section \ref{sec:sloshing} we will report on the sensitivity of the constraint accuracy depending on the accuracy of the conjugate gradient solver and show the benefits of using the proposed time integration method including elimination of constraint drift.

\subsection{Summary}
In this section we have outlined the fourth and main novelty in this paper: a new \textit{constraint-consistent}, high-order accurate time integration method. We have used the \textit{discretize first, substitute next} principle to construct a half-explicit Runge-Kutta method that in a discrete sense possesses the same constraints as the continuous and semi-discrete equations, as derived in sections \ref{sec:governing_equations} and \ref{sec:temporal_properties}. A custom-made three-stage, third order method has been derived such that the differential-algebraic nature of the problem does not lead to order reduction for time-dependent boundary conditions.

\section{Results}
Three test cases are studied in this section to highlight the properties of the proposed Runge-Kutta time integration strategy. The test cases exhibit an increasing level of difficulty in the type of boundary conditions:
\begin{itemize}
\item The growth of Kelvin-Helmholtz instabilities on a periodic domain to study the order of accuracy of the time integration methods.
\item Sloshing of liquid in a closed tank to study solid wall boundary conditions and to assess constraint accuracy and conservation properties in the presence of shock waves.
\item The ramp-up of the gas flow rate in a pipeline to study the order of accuracy for time-dependent boundary conditions, highlighting our proposed RK3 method, which does not suffer from order reduction.
\end{itemize}

\subsection{Kelvin-Helmholtz instability}
The Kelvin-Helmholtz instability occurs due to an imbalance between inertial forces (destabilizing) and gravity forces (in terms of density differences, stabilizing). In multiphase flow in pipelines this instability can start from stratified flow and lead to the formation of slug flow.
We consider a case similar to the one considered by Liao et al.\ \cite{Liao2008}, for conditions where the two-fluid model is \textit{unstable and well-posed}, so we can study the growth of waves that can lead to slug formation. First we find a steady state solution for the parameter values given in table \ref{tab:VKH}. Choosing $u_{l} = 1 \text{ m/s}$ and $\alpha_{l}=0.9$ yields the gas velocity and the pressure gradient necessary to sustain the steady solution:
\begin{equation}
u_{g} = 8.0 \text{ m/s}, \qquad \frac{\rd p_{\text{body}}}{\rd s} = - 87.9 \text{ Pa/m}.
\end{equation}
The velocity difference $u_{g}-u_{l}$ is just below the limit given by equation \eqref{eqn:IKH_limit}, which means that the initial boundary value problem is well-posed. 


\begin{table}[hbtp]
\centering
\caption{Parameter values for the test case with the Kelvin-Helmholtz instability. \label{tab:VKH}}
\begin{tabular}{lrl}
\toprule
parameter & value & unit \\
\midrule
$\rho_{l}$  & $1000$ & \si{kg/m^{3}} \\ 
$\rho_{g}$  & $1.1614$ & \si{kg/m^{3}} \\ 
$R$ & $0.039$ & \si{m} \\
$p_{\text{outlet}}$ & $10^6$ & \si{N/m^{2}} \\ 
$g$ & $9.8$ & \si{m/s^{2}} \\
$\mu_{g}$ & $1.8 \cdot 10^{-5}$ & \si{Pa.s} \\
$\mu_{l}$ & $8.9 \cdot 10^{-4}$ & \si{Pa.s} \\ 
$\epsilon$ & $10^{-8}$ & \si{m} \\
$L$ & $1$ & \si{m} \\
\bottomrule
\end{tabular}
\end{table}

Secondly, we perturb the steady state by imposing a sinusoidal disturbance with wavenumber $k=2\pi$ and a very small amplitude. Linear stability analysis \cite{Liao2008,Sanderse2017} gives the following angular frequencies $\omega$:
\begin{align}\label{eqn:omega_exact}
\omega_{1} & \approx 3.22 +  2.00 i \, \text{s}^{-1}, \\
\omega_{2} &\approx 10.26 - 1.61 i \, \text{s}^{-1}.
\end{align}
The negative imaginary part of $\omega_{2}$ makes the solution unstable due to exponential growth in time. The linear stability analysis provides an analytical solution involving two waves with these frequencies. We set the perturbation related to $\omega_{1}$ to zero, implying that a single wave with frequency $\omega_{2}$ results. The exact solution to the linearized equations is then
\begin{equation}\label{eqn:analytical_linearized}
\vt{W} (s,t) = \vt{W}^{0} + \operatorname{Re} \left[ \vt{\varepsilon}_2 e^{i(\omega_2 t - ks)} \right],
\end{equation}
where  $\vt{\varepsilon}_2$ is obtained by choosing the liquid hold-up fraction perturbation as $\hat{\alpha}_{l}=10^{-6}$, and then computing the perturbations in the gas and liquid velocity from the dispersion analysis \cite{Liao2008}. The initial condition which follows by taking $t=0$ does not satisfy \eqref{eqn:IC2}-\eqref{eqn:IC3} exactly, and we therefore perform a projection step to make the initial conditions consistent.

%
As a consistency check, we have first investigated the accuracy of the entire space-time discretization by comparing the discrete solution to the linearized analytical solution \eqref{eqn:analytical_linearized} at $t=1$. Upon refining simultaneously the grid and the time step, second order accuracy is obtained for all solution components, for RK2, the proposed RK3 method, and RK4. This is because the convergence rate is dominated by the second order accuracy of the spatial discretization.

\begin{figure}[hbtp]
\centering
	\begin{subfigure}[b]{.49\textwidth}
	\centering
	\includegraphics[width=\textwidth]{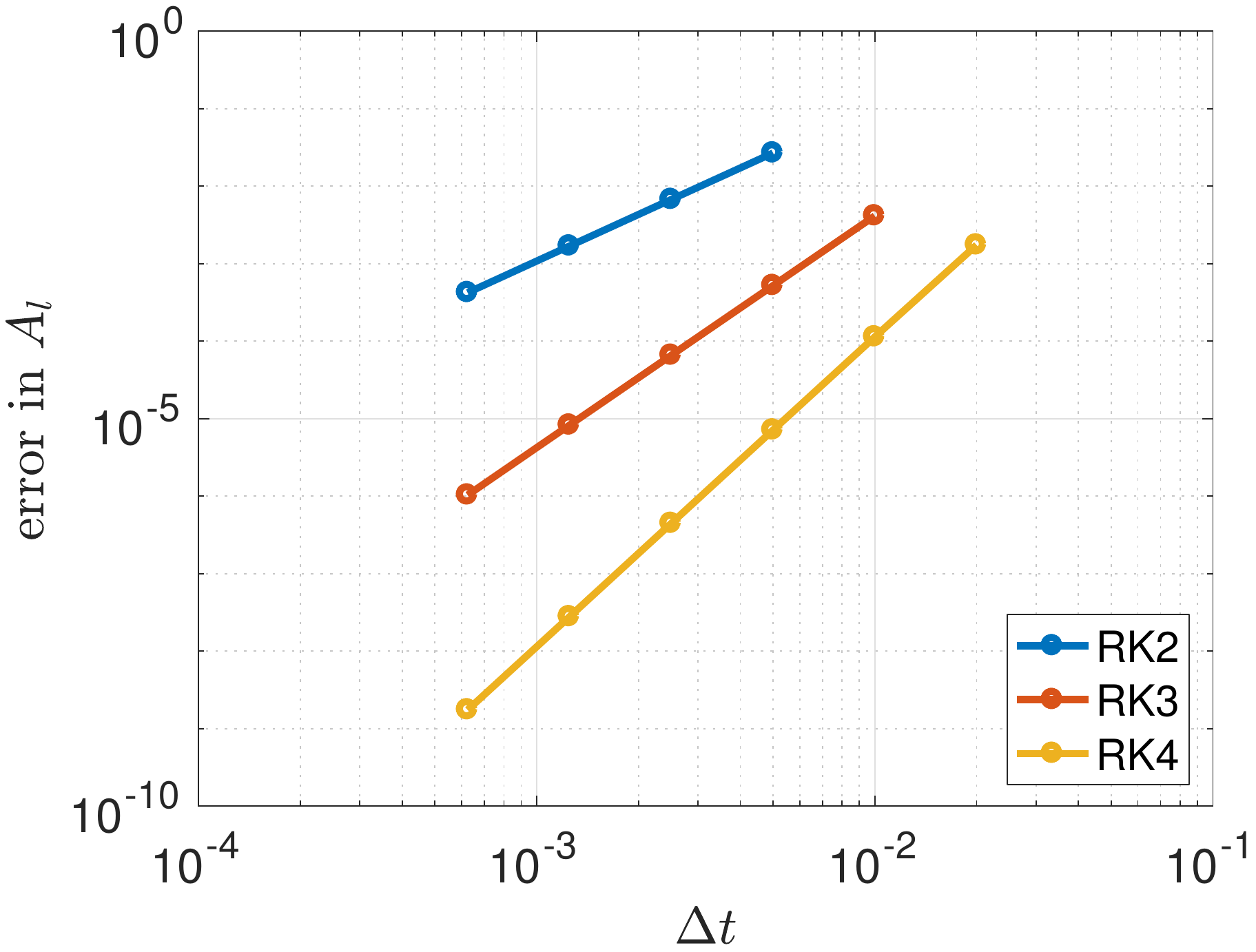}
	\caption{Hold-up fraction.}
	\end{subfigure}
	\hfill
	\begin{subfigure}[b]{.49\textwidth}	
	\centering
	\includegraphics[width=\textwidth]{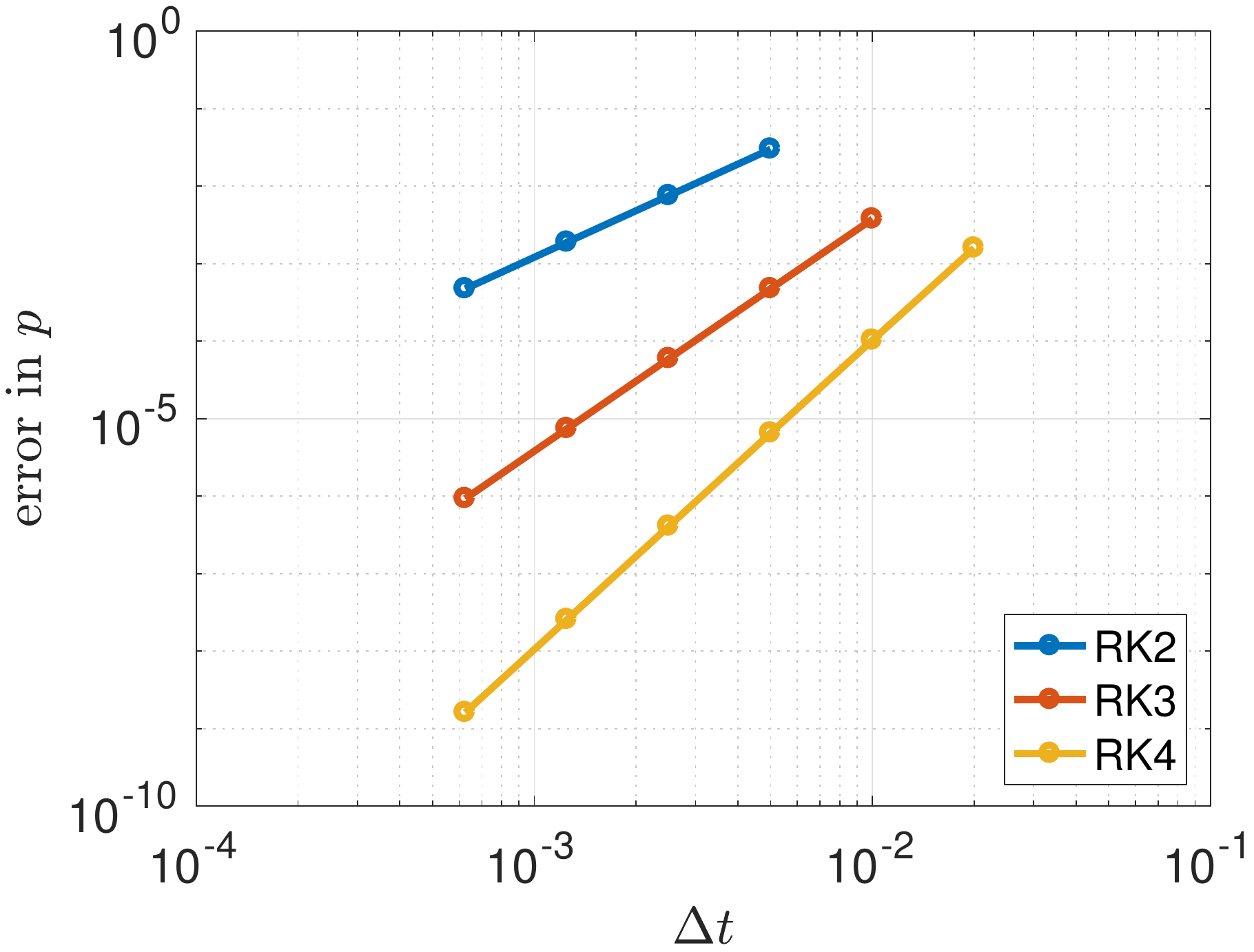}
	\caption{Pressure.}
	\end{subfigure}\\
	\vspace{0.5cm}
	\begin{subfigure}[b]{.49\textwidth}
	\centering
	\includegraphics[width=\textwidth]{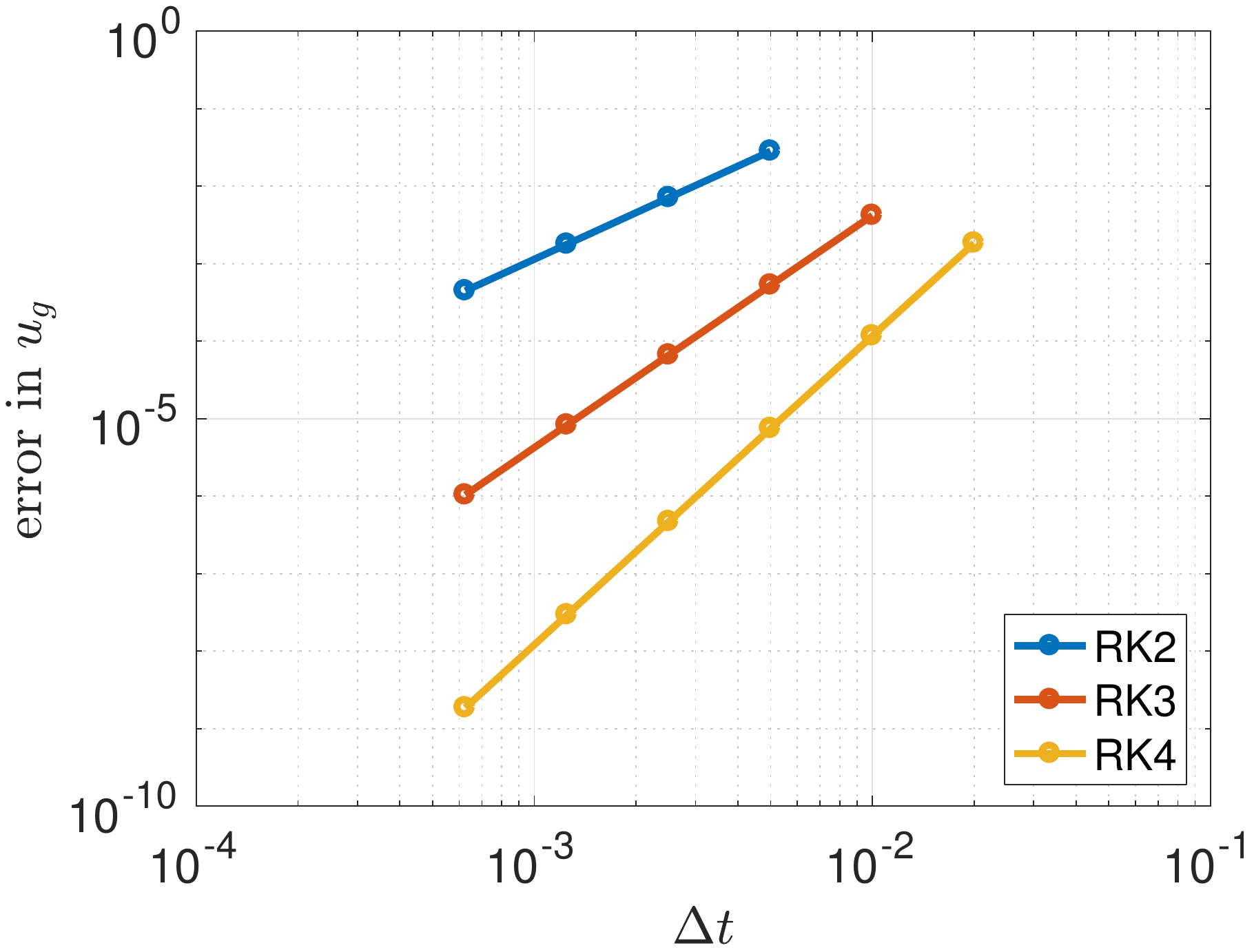}
	\caption{Gas velocity.}
	\end{subfigure}
	\hfill
	\begin{subfigure}[b]{.49\textwidth}	
	\centering
	\includegraphics[width=\textwidth]{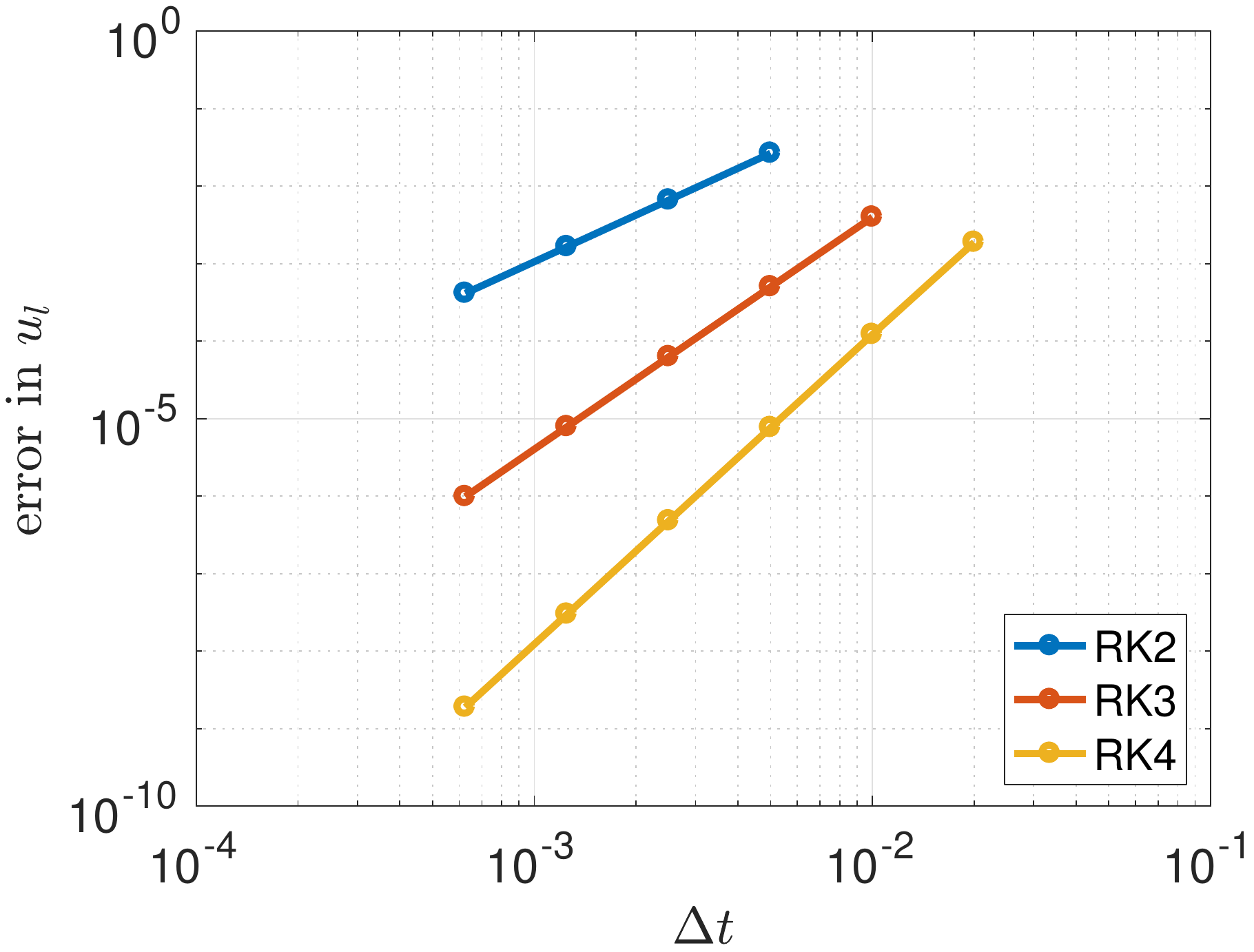}
	\caption{Liquid velocity.}
	\end{subfigure}
\caption{Convergence of the temporal error for the incompressible Kelvin-Helmholtz problem.\label{fig:KH_INC_errors}}
\end{figure}

To investigate the temporal accuracy alone, we compute a reference solution $\vt{W}_{\text{ref}}$ at $t=1$ \si{s} with $N=40$ volumes and RK4 with a small time step, $\Delta t= 1 \cdot 10^{-4}$ \si{s}, so that the temporal error is negligible. The temporal error in the hold-up fraction (scaled by the magnitude of the disturbance) then follows from
\begin{equation}\label{eqn:temporal_error}
\epsilon_{A_{l}, \Delta t} = \frac{ \| A_{l,\Delta t} - A_{l,\text{ref}} \|_{\infty}}{ \hat{A}_{l}},
\end{equation}
and similarly for the other solution components. The perturbation is increased to $\hat{\alpha}_{l}=10^{-3}$ to avoid errors that are around machine precision. Figure \ref{fig:KH_INC_errors} shows these errors as a function of the timestep. All methods convergence to their design order of accuracy, for both the differential and algebraic components. In this test case there are no boundary conditions (so $r=0$) and therefore \textit{no} additional order conditions appear, making RK4 an excellent choice.

\FloatBarrier

\subsection{Sloshing in a closed cylindrical tank}\label{sec:sloshing}
In this section we present a new and challenging test case for the incompressible two-fluid model: sloshing of liquid and gas in a closed pipe section. In a closed system the mass of both the gas and liquid phase is conserved exactly, and we desire the same property of our numerical algorithm. Furthermore, this test case is well-suited for checking the accuracy of the volume constraint C0 and volumetric flow constraint C1, because any constraint drift or error growth will not be able to `leave' the domain since there are no outflow boundaries.

The geometry of the problem is shown in figure \ref{fig:sloshing_tank}. At $t=0$, the liquid is released and starts flowing towards the left wall, from which it reflects, resulting in a complicated wave pattern in time and space. In an experimental setting this can be accomplished by suddenly tilting the pipe section from its horizontal position, as is done for example in \cite{Thorpe1969} for a liquid-liquid system. The parameter values that we employ are shown in table \ref{tab:closedtank_parameters}. The spatial mesh has $N=80$ finite volumes.

\begingroup
\begin{figure}[hbtp]
\fontfamily{lmss}
\fontsize{10pt}{12pt}\selectfont
\centering 
\def\svgwidth{0.7 \textwidth}
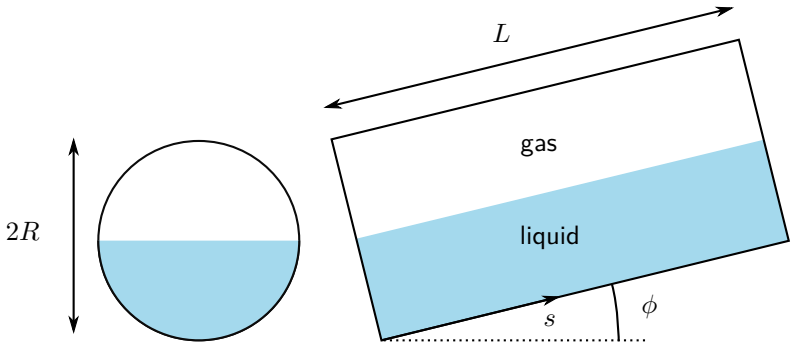 
\caption{Initial condition for sloshing simulation.}
\label{fig:sloshing_tank}
\end{figure}
\endgroup

\begin{table}[hbtp]
\centering
\caption{Parameter values for the closed tank problem. \label{tab:closedtank_parameters}}
\begin{tabular}{lrl}
\toprule
parameter & value & unit \\
\midrule
$\rho_{l}$  & $1000$ & \si{kg/m^{3}} \\ 
$\rho_{g}$  & $1.1614$ & \si{kg/m^{3}} \\ 
$R$ & $0.05$ & \si{m} \\
$g$ & $9.8$ & \si{m/s^{2}} \\
$\mu_{g}$ & $1.5 \cdot 10^{-2}$ & \si{Pa.s} \\
$\mu_{l}$ & $5.0 \cdot 10^{-2}$ & \si{Pa.s} \\ 
$\epsilon$ & $10^{-8}$ & \si{m} \\
$L$ & $1$ & \si{m} \\
$\phi$ & $2$ & \si{deg} \\
\bottomrule
\end{tabular}
\end{table}

\subsubsection{Qualitative analysis}
\begin{figure}[hbtp]
\centering
	\includegraphics[width=\textwidth]{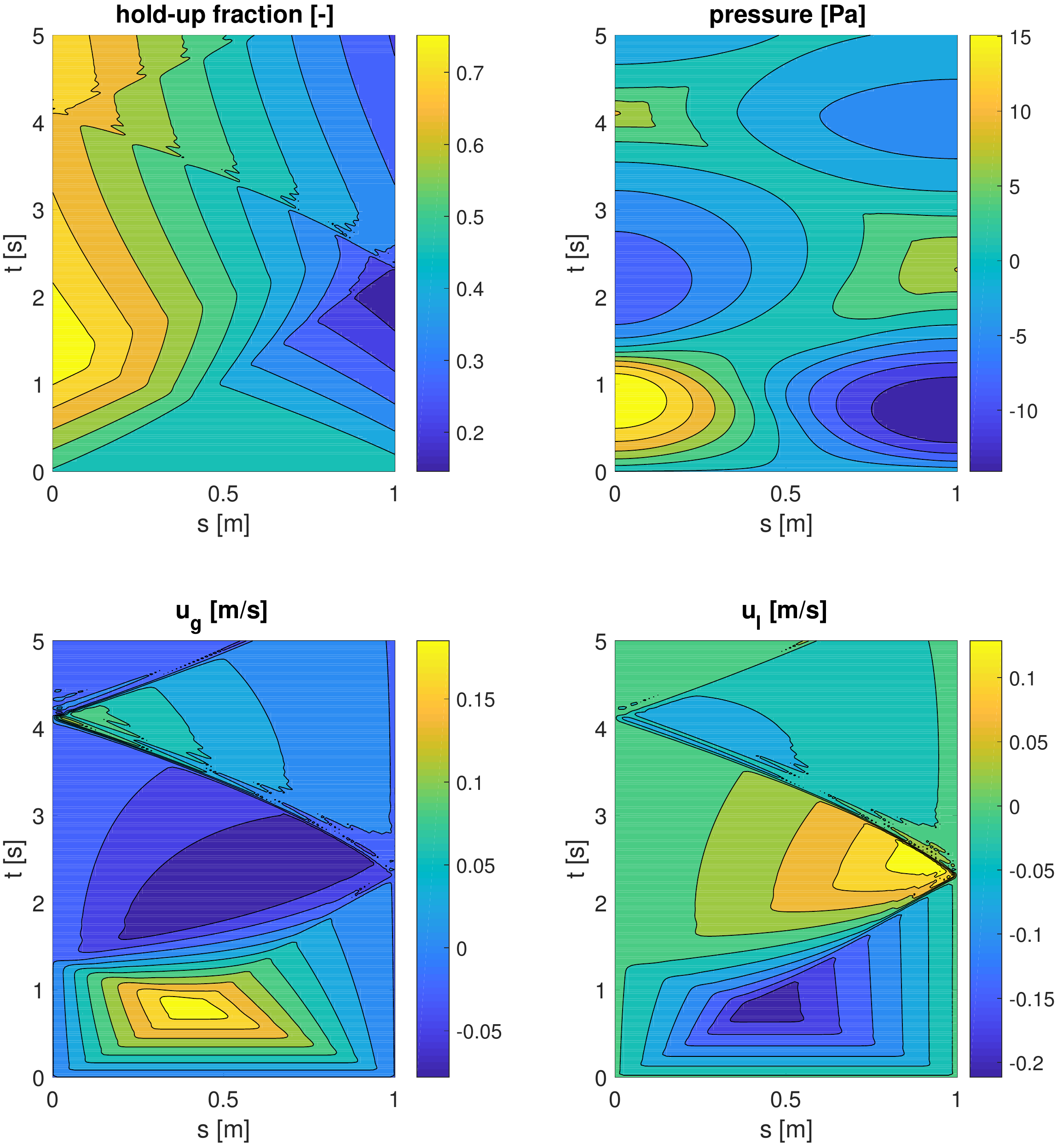}
\caption{Space-time solution for the sloshing problem with $N=80$, $\Delta t=0.02$ and RK4.\label{fig:closed_tank_fullsolution}}
\end{figure}

\begin{figure}[hbtp]
\centering
	\includegraphics[width= \textwidth]{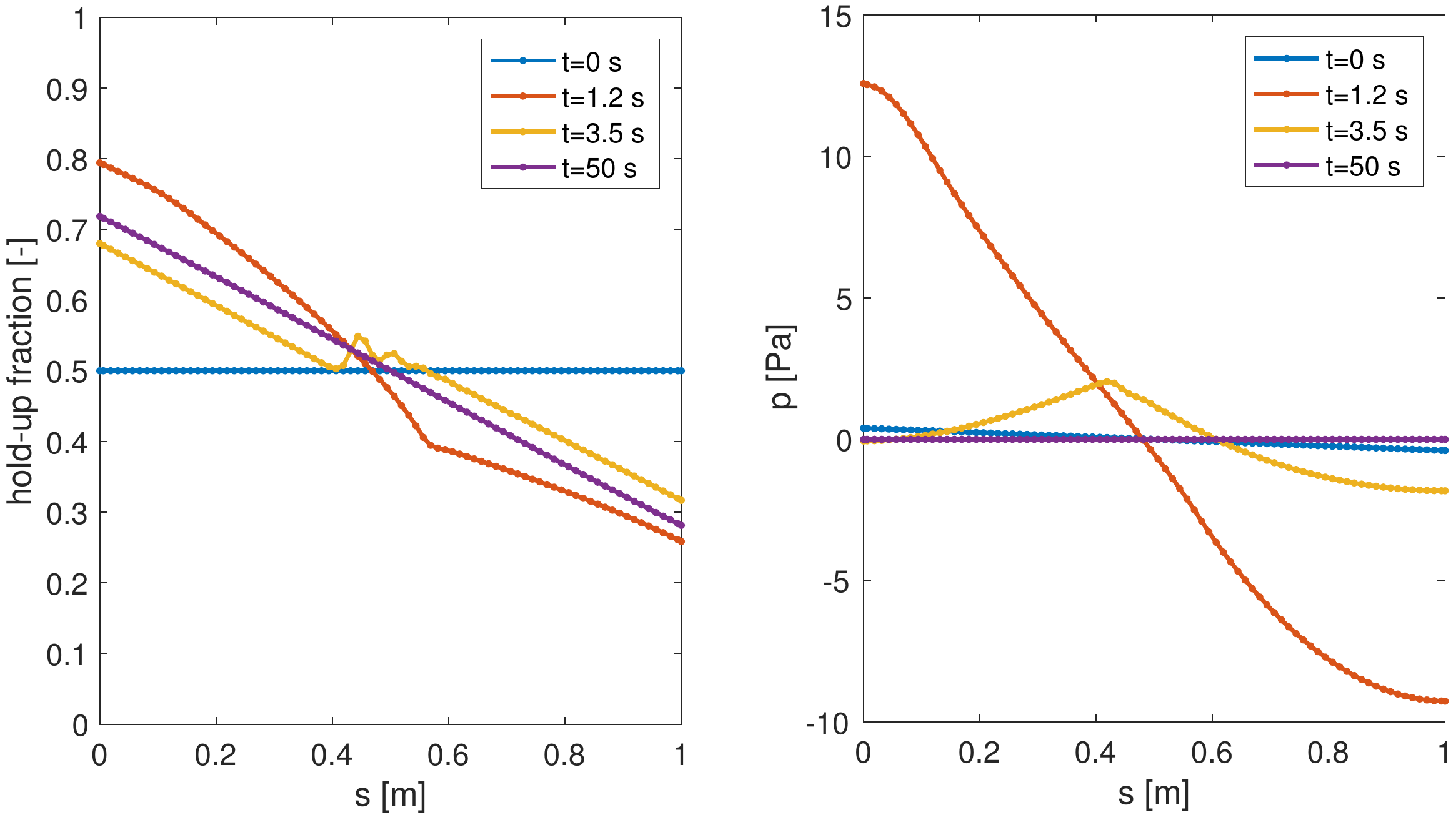}
\caption{Details from figure \ref{fig:closed_tank_fullsolution}: initial, intermediate (with shock wave at $t=3.5$ \si{s}) and final solutions for the hold-up fraction (left) and the pressure (right).\label{fig:closed_tank_steadystate}}
\end{figure}

Figure \ref{fig:closed_tank_fullsolution} shows the solution in terms of hold-up fraction, pressure and phase velocities for the first 5 seconds of the simulation. Snapshots of the hold-up fraction and pressure are given in figure \ref{fig:closed_tank_steadystate}. 

A uniform initial condition for all parameters satisfies the volume constraint and volume flow constraint, but not the pressure constraint. The initial condition for the pressure should be determined from equation \eqref{eqn:IC3}, which results in a pressure that changes linearly along the $s$ direction, due to the effect of gravity (see figure \ref{fig:closed_tank_steadystate} right, blue line). Note that the pressure behaviour at the boundary is correct without requiring an explicit pressure boundary condition.

The sudden tilting of the pipe section creates two waves, originating from the two boundaries. A compression wave moves to the right, increasing the hold-up fraction and pressure on the left side of the tank. An expansion wave moves to the left, decreasing the pressure and hold-up fraction on the right side of the tank. The expansion wave moves slightly faster than the compression wave and reaches the left side after approximately 1.3 seconds, where it reflects. The compression wave reaches the right side after approximately 2.4 seconds, continues to steepen, and forms a shock wave (this is particularly clear from the plots of the hold-up fraction and the gas velocity). Since central differences are used for the spatial discretization, small wiggles in the solution are present close to the shock front. \hl{Other spatial discretization methods, e.g.\ a Roe-scheme \cite{Akselsen2017}, could be employed in order to resolve shocks without oscillations. In principle such methods can be directly used with our time-integration method, as long as the discrete coupling between mass and momentum equations is satisfied, as mentioned in section \ref{sec:spatial_discretization}}. A detail of the solution at $t=3.5$ seconds is shown in figure \ref{fig:closed_tank_steadystate}, highlighting the shock wave in the hold-up fraction around $s=0.5$ \si{m}. The pressure, on the other hand, does not exhibit a jump in the solution, but contains a jump in its first derivative, due to the elliptic nature of the Poisson equation.

Without further agitation, the sloshing liquid comes to rest due the action of friction. In figure \ref{fig:closed_tank_steadystate} the final steady state (obtained at $t=50$, when $|u_{g}|$, $|u_{l}| < 10^{-8}$) is displayed. In contrast to the initial conditions, the final steady state condition has a uniform pressure value, even though there is more liquid (and therefore hydrostatic head) at the left side of the tank than at the right side. This is because the pressure in the two-fluid model is the \textit{pressure at the interface}. Evaluating the pressure gradient equation, equation \eqref{eqn:constraint3}, for the final quiescent steady state, gives
\begin{equation}
\left( \frac{A_{g}}{\rho_{g}} + \frac{A_{l}}{\rho_{l}} \right) \dd{p}{s} = -A g_{n} \dd{h}{s}  - A g_{s}.
\end{equation}
Since the steady state attains a level surface, the slope of the liquid is $\dd{h}{s} = - \tan \phi$, and the equation above reduces to 
\begin{equation}
\left( \frac{A_{g}}{\rho_{g}} + \frac{A_{l}}{\rho_{l}} \right) \dd{p}{s} = 0.
\end{equation}
In physical terms, the hydrostatic head of the liquid is counterbalanced by the level gradient term. Of course, this holds only at the interface; in the lower left corner of the pipe section, for example, the pressure will be higher. 

\FloatBarrier

\subsubsection{Conservation properties and order of accuracy}
\begin{figure}[hbtp]
\centering
	\begin{subfigure}[b]{.49\textwidth}
	\includegraphics[width=\textwidth]{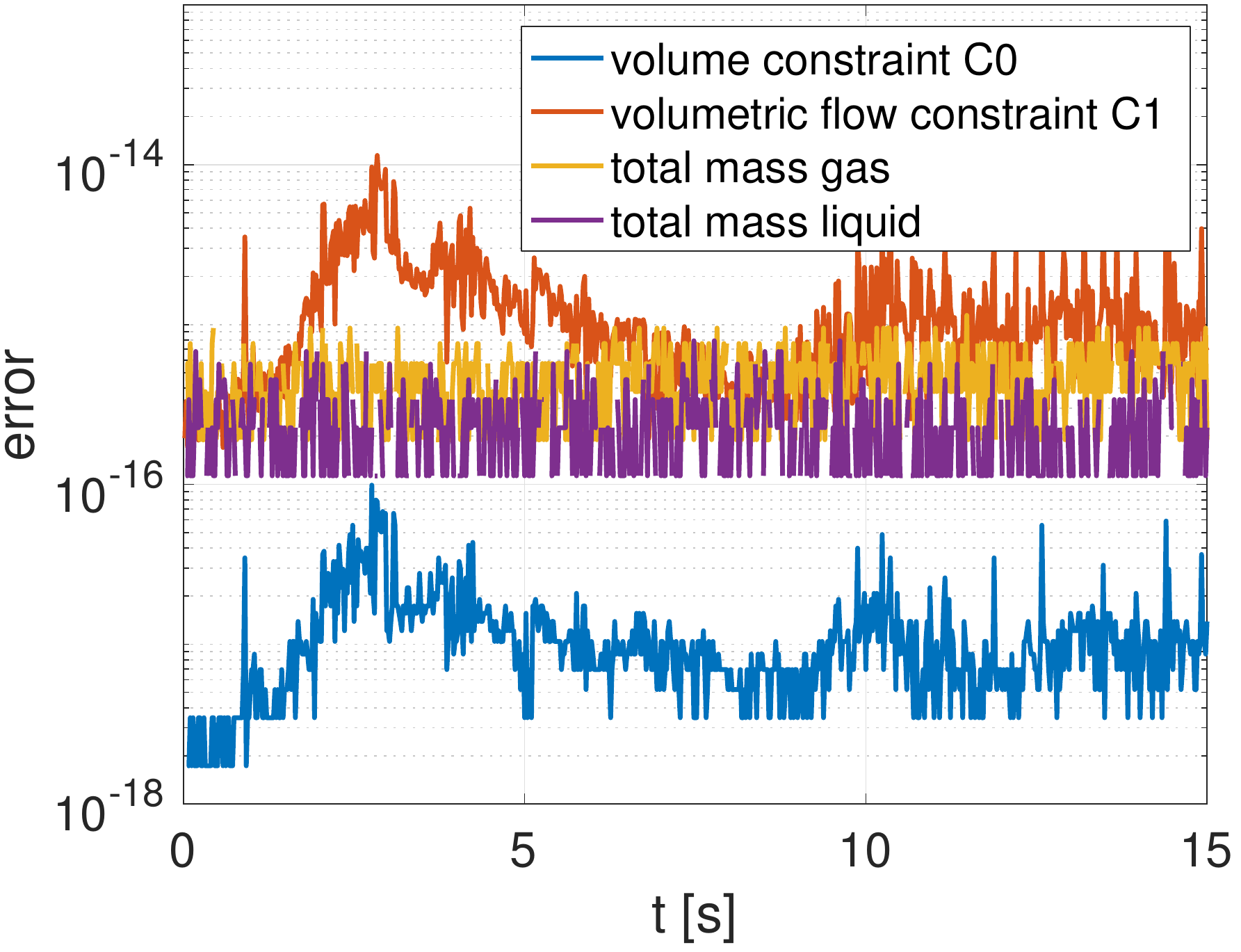}
	\caption{Tolerance $10^{-12}$. \label{fig:closed_tank_conservation_a}}
	\end{subfigure}
	\begin{subfigure}[b]{.49\textwidth}
	\centering
	\includegraphics[width=\textwidth]{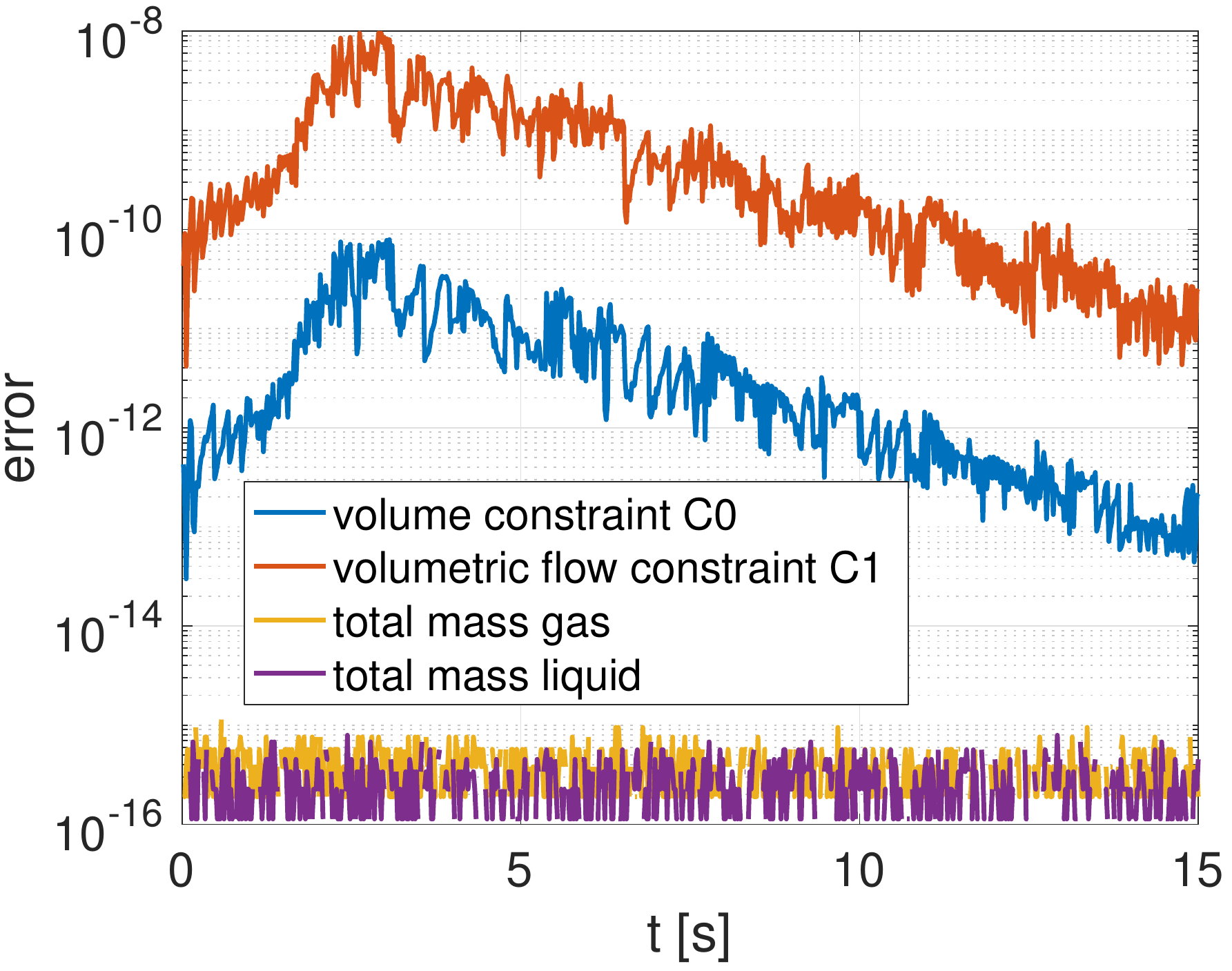}
	\caption{Tolerance $10^{-6}$. \label{fig:closed_tank_conservation_b}}
	\end{subfigure}\\
	\vspace{0.5cm}
	\begin{subfigure}[b]{.49\textwidth}
	\centering
	\includegraphics[width=\textwidth]{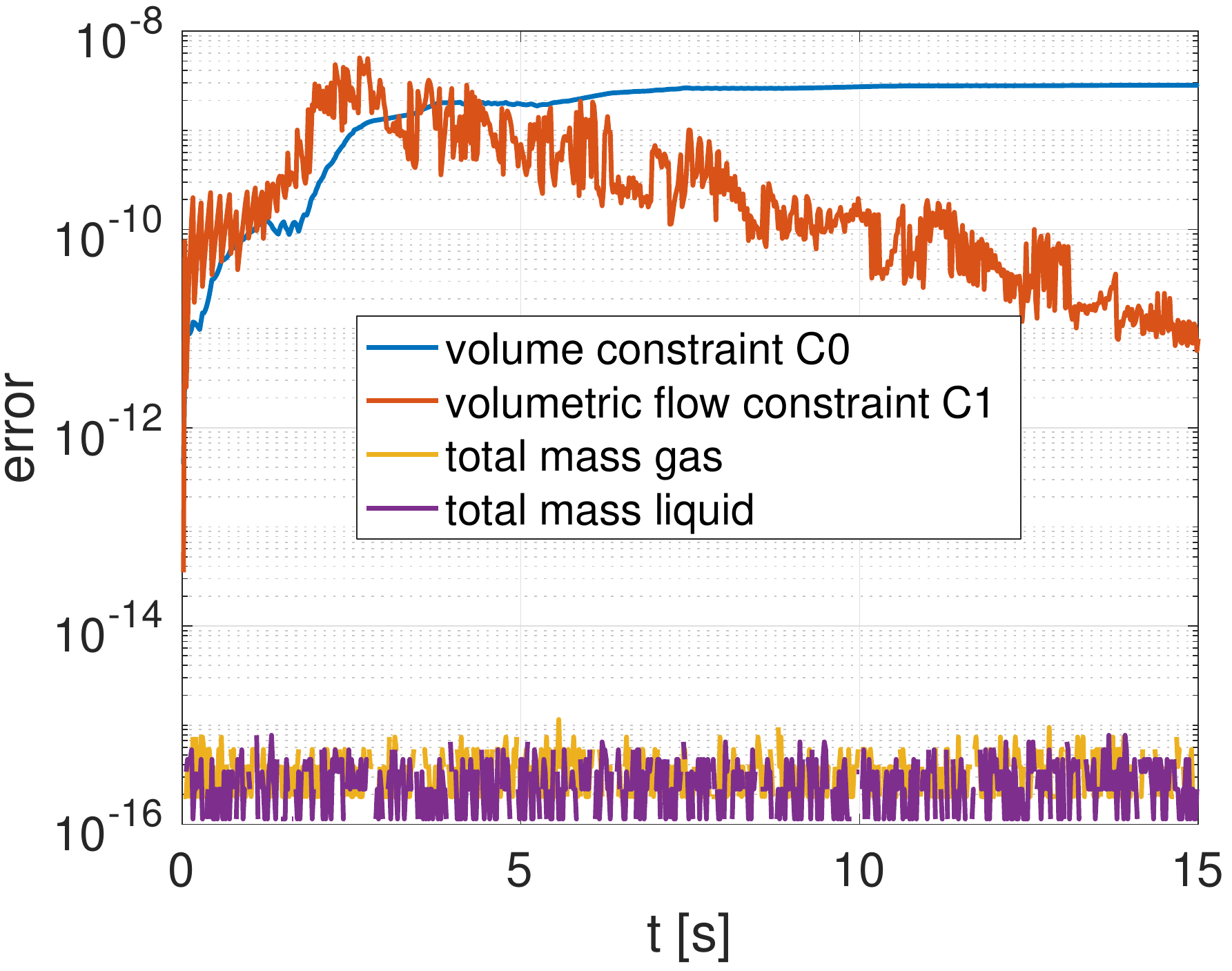}
	\caption{Tolerance $10^{-6}$, without correction terms. \label{fig:closed_tank_conservation_c}}
	\end{subfigure}
\caption{Error in conservation and constraint properties depending on the tolerance of the  conjugate gradient solver, for the sloshing problem with $N=80$, $\Delta t=0.02$ and RK4.\label{fig:closed_tank_conservation}}
\end{figure}

Figure \ref{fig:closed_tank_conservation_a} highlights the conservation and constraint properties of our time integration method. When the tolerance of the conjugate gradient solver is sufficiently small, both the volume constraint C0 and the volumetric flow constraint C1 are satisfied until machine precision over the entire course of the simulation. When the threshold of the conjugate gradient solver is increased, the errors in the constraint are larger, but do not drift due to the correction terms proposed in equation \eqref{eqn:ppe_correction} (see figure \ref{fig:closed_tank_conservation_b}). In contrast, when the correction terms as proposed in equation \eqref{eqn:poisson_additional} are not taken into account, the volume constraint starts drifting and the drift does not decrease when the solution approaches a steady state (figure \ref{fig:closed_tank_conservation_c}). Independent of conjugate gradient solver tolerance and correction terms in the pressure equation, the mass of both the gas phase and the liquid phase is conserved until machine precision in all cases (this is a property of the finite volume method).

\begin{figure}[hbtp]
\centering
	\begin{subfigure}[b]{.49\textwidth}
	\centering
	\includegraphics[width=\textwidth]{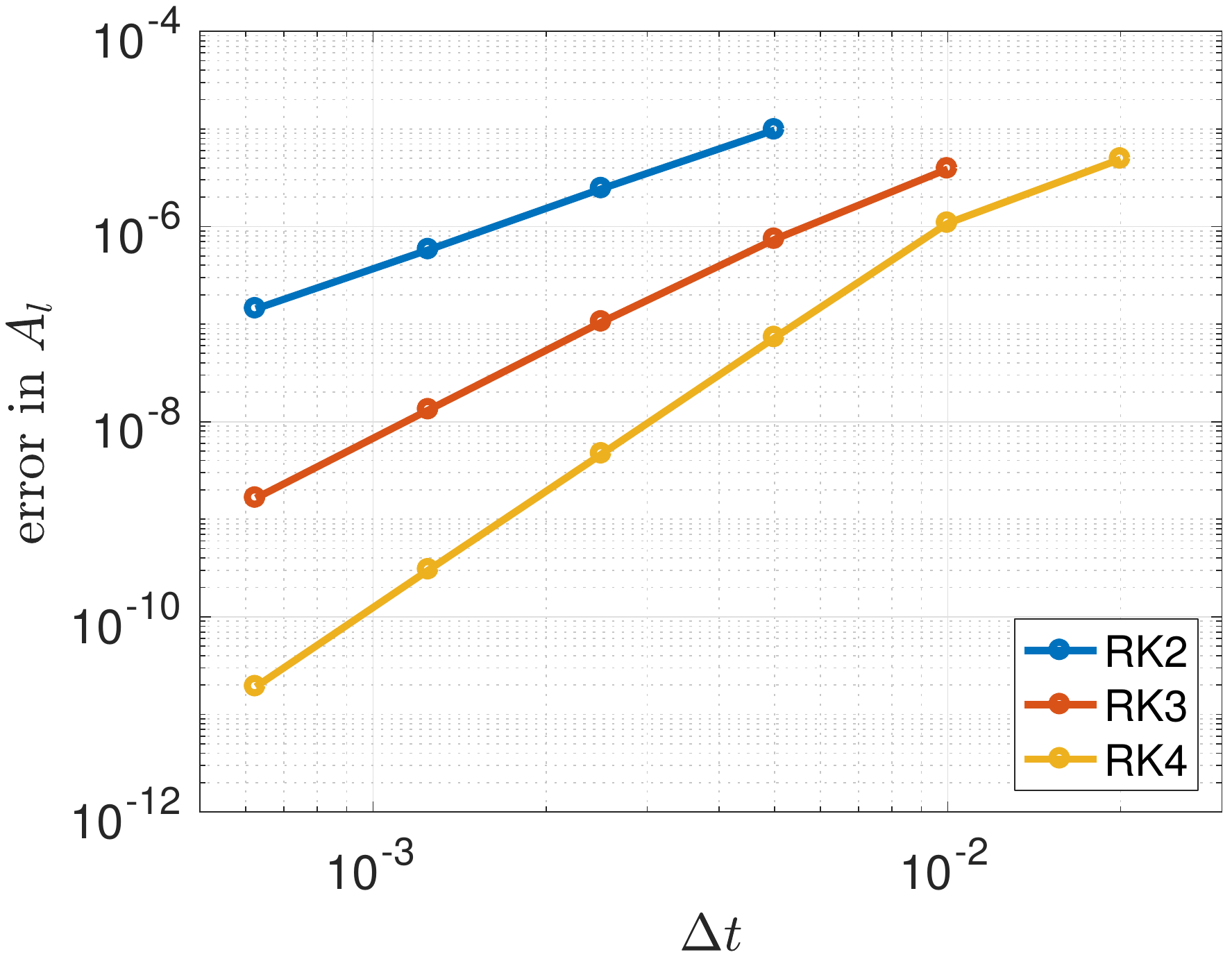}
	\caption{Hold-up fraction.}
	\end{subfigure}
	\hfill
	\begin{subfigure}[b]{.49\textwidth}	
	\centering
	\includegraphics[width=\textwidth]{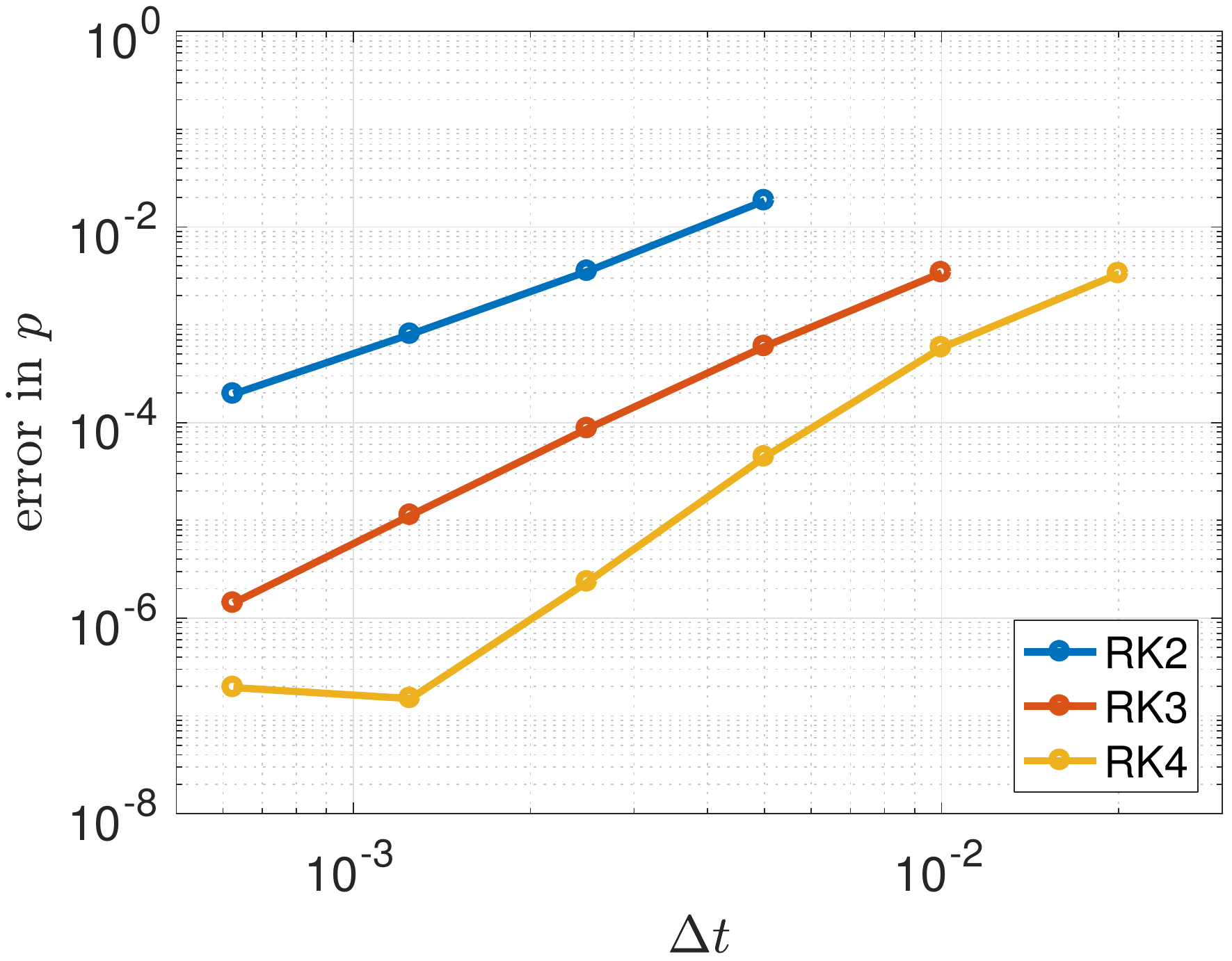}
	\caption{Pressure.}
	\end{subfigure}\\
	\vspace{0.5cm}
	\begin{subfigure}[b]{.49\textwidth}
	\centering
	\includegraphics[width=\textwidth]{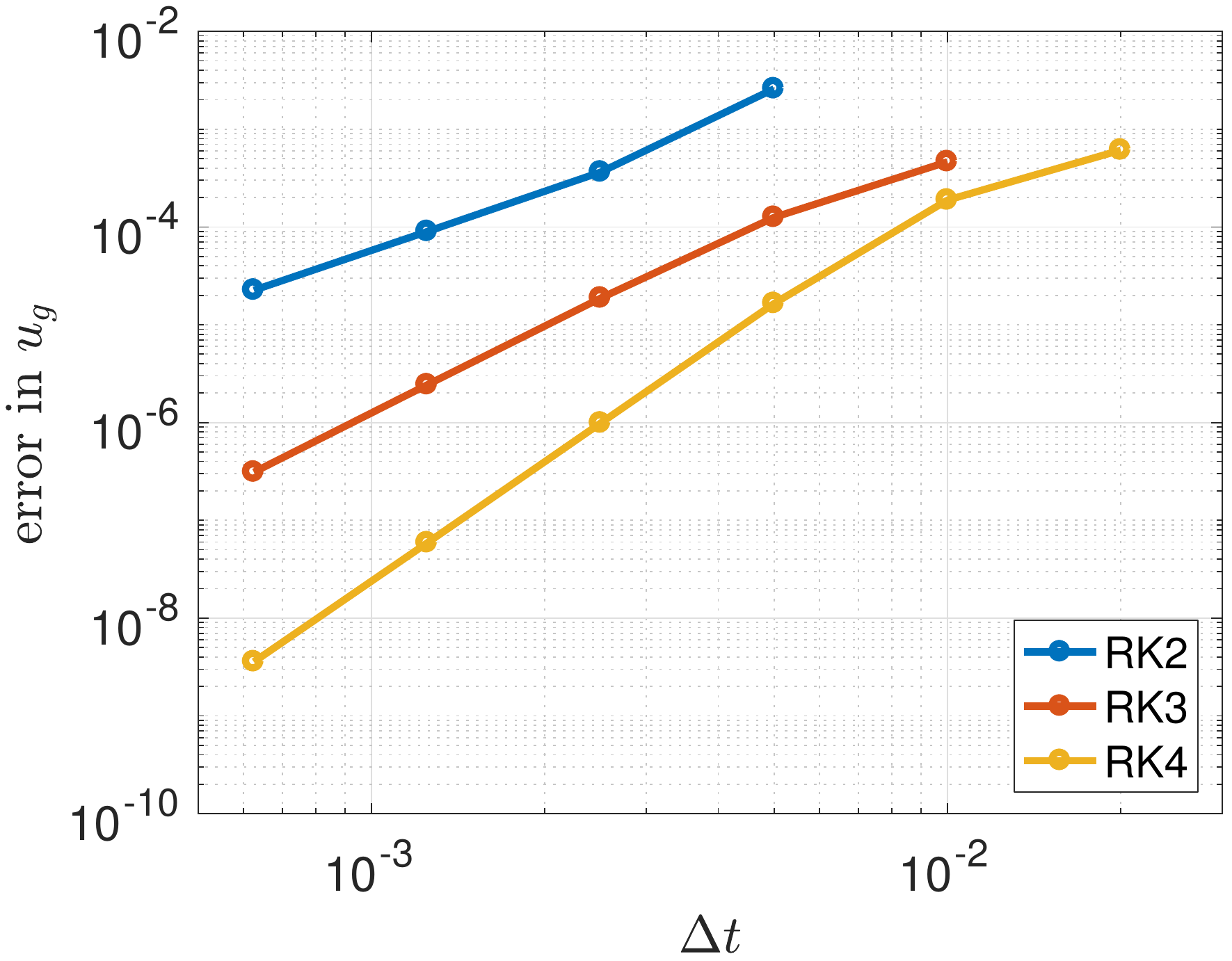}
	\caption{Gas velocity.}
	\end{subfigure}
	\hfill
	\begin{subfigure}[b]{.49\textwidth}	
	\centering
	\includegraphics[width=\textwidth]{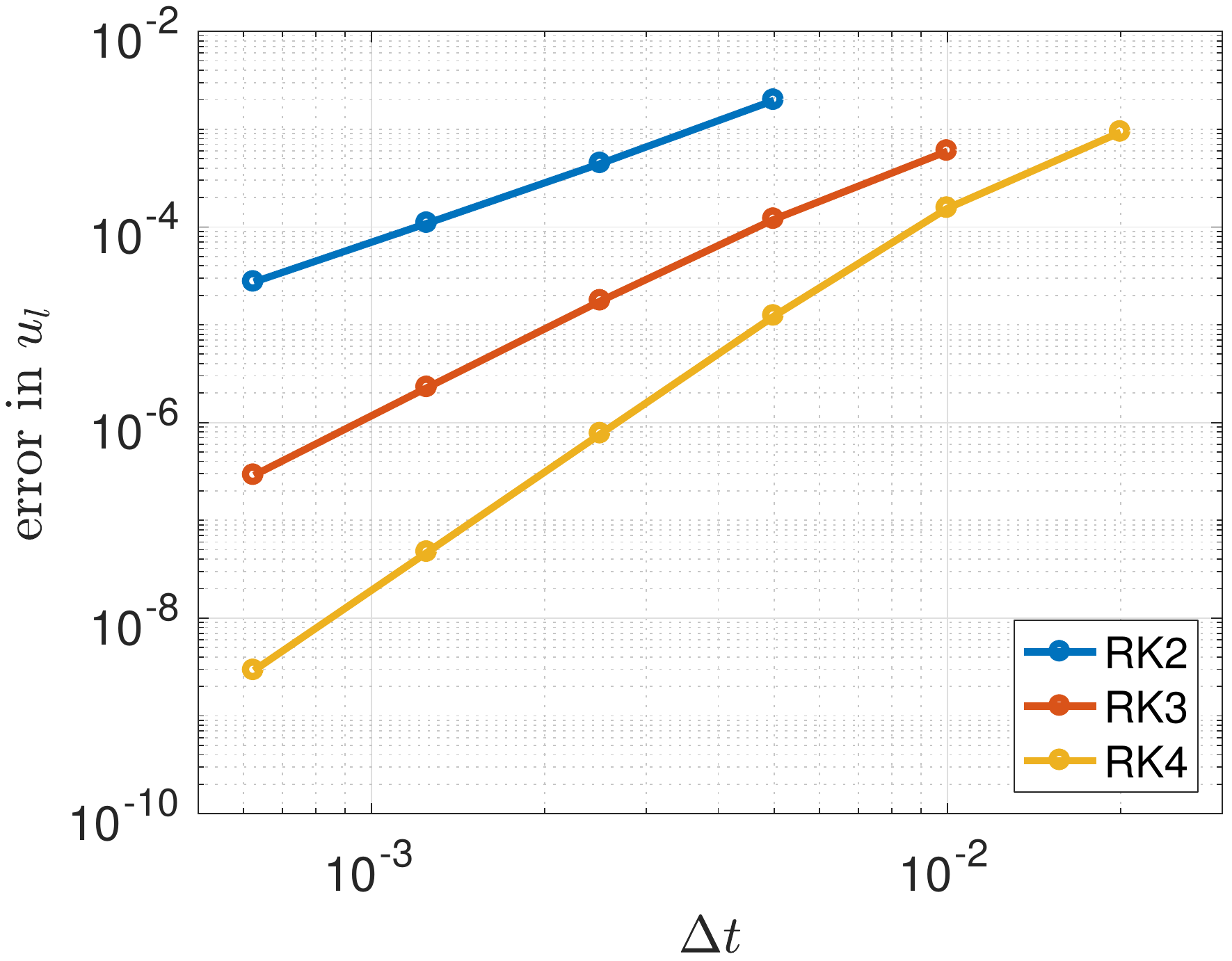}
	\caption{Liquid velocity.}
	\end{subfigure}
\caption{Convergence of the temporal error for the sloshing problem.\label{fig:closed_tank_convergence_study}}
\end{figure}

A temporal accuracy study is performed in the same way as for the Kelvin-Helmholtz test case. We use a reference solution obtained with RK4 and $\Delta t=10^{-4}$ \si{s} to compute the temporal error at $t=1.2 \si{s}$. The resulting convergence of the liquid hold-up fraction, phase velocities, and pressure is shown in figure \ref{fig:closed_tank_convergence_study}. RK2, RK3 and RK4 all converge according to their design order of accuracy. High order accuracy for the pressure is obtained via the solution of equation \eqref{eqn:poisson_additional}. Besides the high accuracy offered by RK4, we observe here another main advantage: RK4 can be used at larger time steps than RK2 and RK3, due to its larger stability domain. Similar to the first test case, also this test case features $r=0$, because the boundary conditions for $I_{g}$ and $I_{l}$ are independent of time. Therefore, there is no order reduction; RK4 achieves its classical order of accuracy and forms an excellent choice for this test case.

\FloatBarrier

\subsection{Perturbed hold-up wave propagation}\label{sec:IFP}

\subsubsection{Test case description}
The last test case we perform is the propagation of a hold-up wave, caused by the increase of the gas production at the inlet of a 1 km multiphase pipeline. The test case is inspired by the one proposed by the French Petroleum Institute (IFP) as described by Omgba-Essama \cite{Omgba-Essama2004}. The parameters of the problem are described in table \ref{tab:IFP_parameters}. In contrast to \cite{Omgba-Essama2004}, we employ the two-fluid model instead of the homogeneous equilibrium mixture (HEM) model, and therefore require different initial conditions in order to prevent ill-posedness \cite{Omgba-Essama2004} (the HEM model is unconditionally well-posed). The initial conditions are steady state production with inlet mass flows of liquid and gas of $I_{l}=1$ \si{kg/s} and $I_{g,\text{start}}=0.02$ \si{kg/s}. Furthermore, instead of a linear change in the gas mass flow rate to $I_{g,\text{end}}=0.04$ \si{kg/s} proposed by \cite{Omgba-Essama2004}, we employ a sinusoidally varying flow rate, smoothly started from the initial conditions:
\begin{equation}
I_{g} = I_{g,\text{start}}  + ( I_{g,\text{end}}- I_{g,\text{start}} ) e^{1-10/t} \left( \frac{1}{2}+\sin(t/5)^2 \right)/e^{1}.
\end{equation}
The period of oscillation is $5 \pi \approx 15.7 \text{s}$. This sine-type inflow provides a more severe testcase than a linear ramp-up, because the term $\dot{r}$ in equation \eqref{eqn:poisson_additional} contains time-dependent terms. 

\begin{figure}[hbtp]
\centering
	\includegraphics[width=\textwidth]{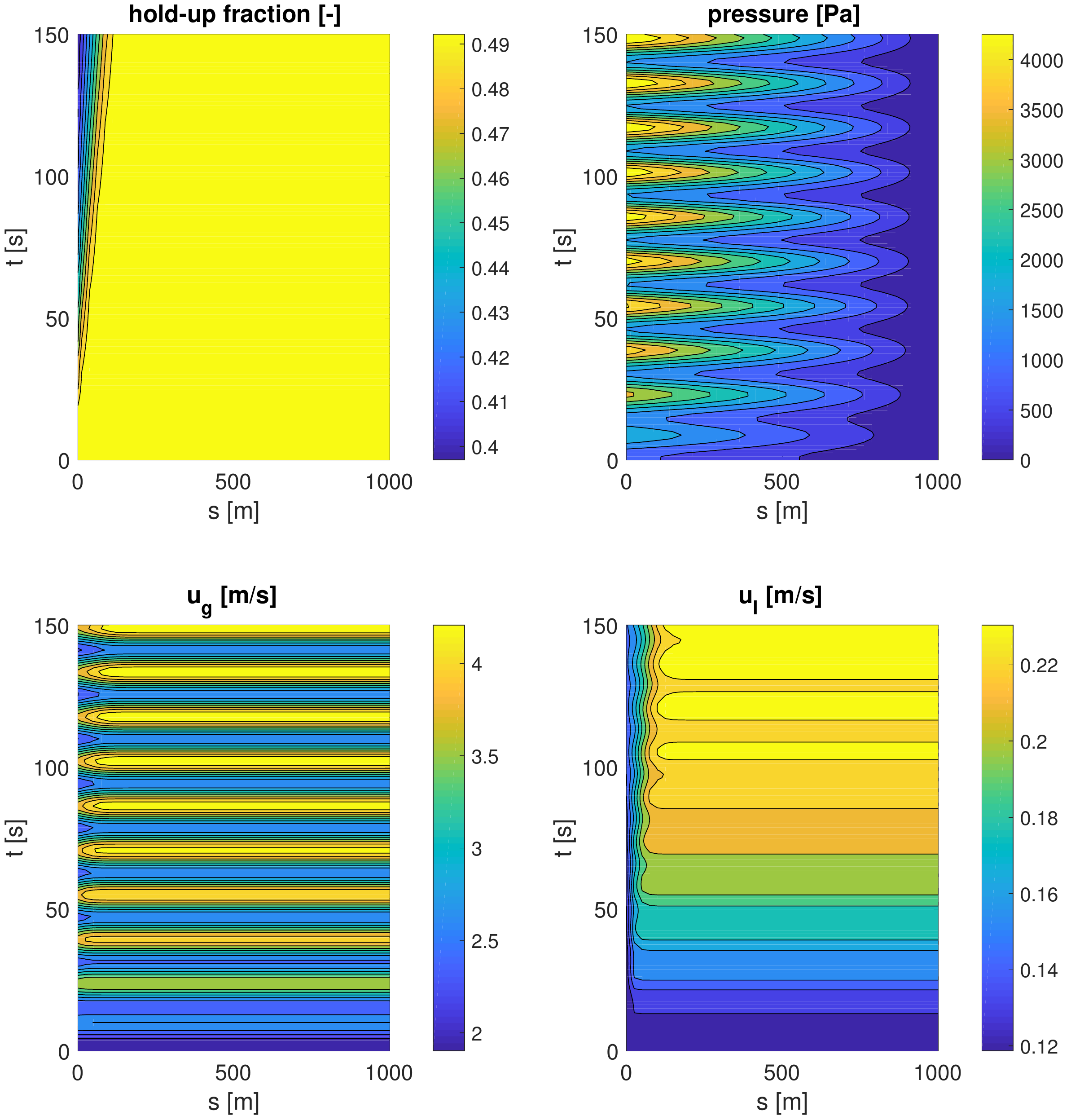}
\caption{Space-time solution for the IFP problem with $N=40$, $\Delta t=1.25$ \si{s} and the proposed RK3 method.\label{fig:IFP_fullsolution}}
\end{figure}

\begin{figure}[hbtp]
\centering
	\includegraphics[width=\textwidth]{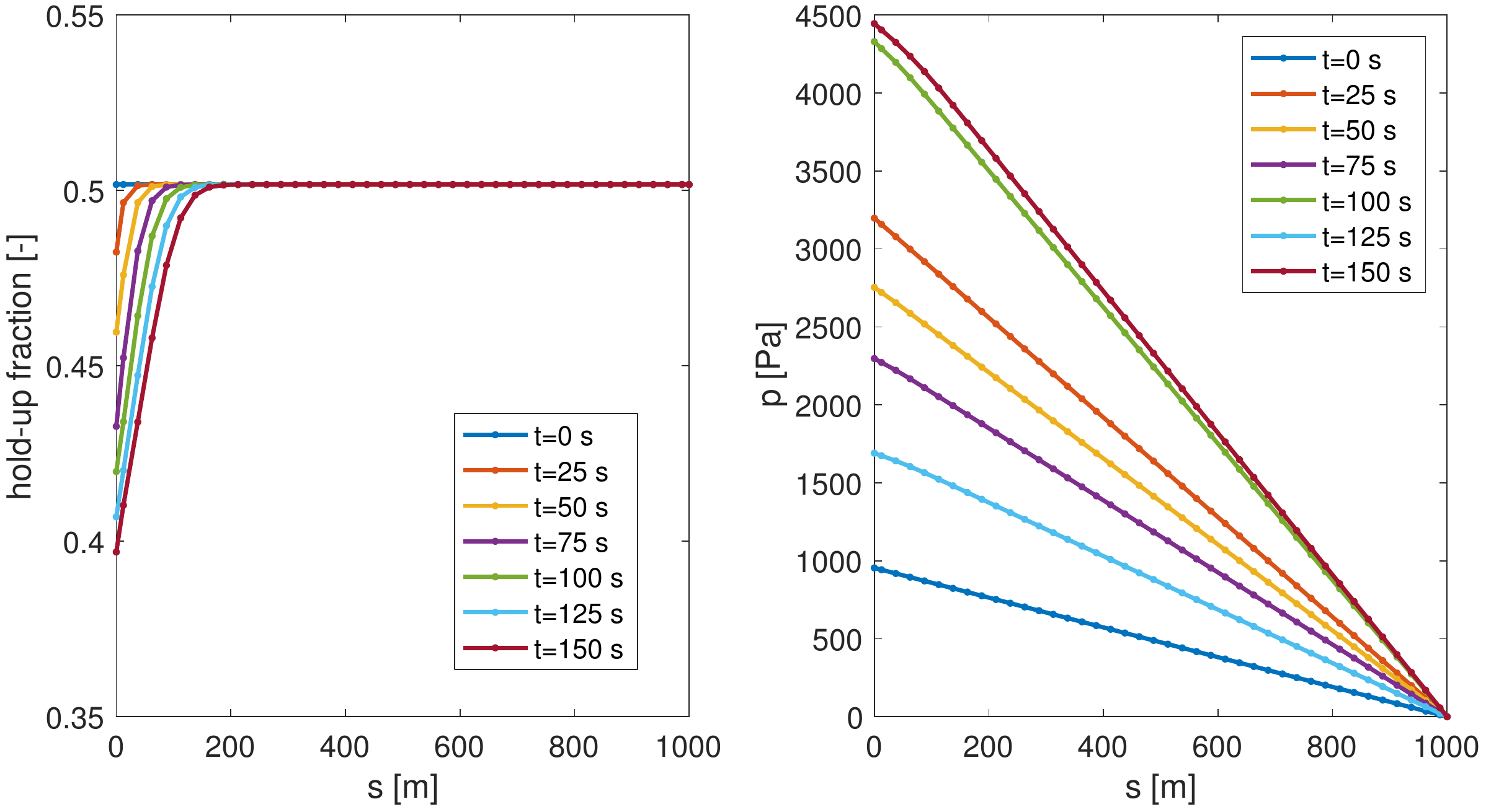}
\caption{\hl{Details of figure \ref{fig:IFP_fullsolution} for specific time instances.\label{fig:IFP_intermediates}}}
\end{figure}

A qualitative view of the solution behaviour in space and time until $t=150$ \si{s} is shown in figures \ref{fig:IFP_fullsolution} \hl{and \ref{fig:IFP_intermediates}}. Since the absolute value of the pressure is not important in incompressible calculations, the pressure difference with respect to $p_{\text{outlet}}$ is shown. At $t=0$, all quantities are uniform in space, except the pressure, which decreases as a function of $s$ due to friction losses. Like in the previous test case, the initial condition for the pressure is not prescribed, but is determined by solving equation \eqref{eqn:IC3}. After a few seconds, the increasing gas mass flow rate leads to a hold-up wave propagating through the pipeline. Two transient effects play a role. First, the gas velocity increases almost instantaneously to adjust for the higher mass flow rate, the liquid velocity increases due to interfacial friction, and the pressure drop increases due to the higher wall friction.  This process is repeated given the periodic nature of the inflow. Second, the hold-up fraction starts to slowly adjust following the convection-type equation \eqref{eqn:conservation_mass}, with a convective velocity determined by the magnitude of the eigenvalues. 


\begin{table}[hbtp]
\centering
\caption{Parameter values for the IFP problem. \label{tab:IFP_parameters}}
\begin{tabular}{lrl}
\toprule
parameter & value & unit \\
\midrule
$\rho_{l}$  & $1003$ & \si{kg/m^{3}} \\ 
$\rho_{g}$  & $1.26$ & \si{kg/m^{3}} \\ 
$R$ & $0.073$ & \si{m} \\
$p_{\text{outlet}}$ & $10^6$ & \si{N/m^{2}} \\ 
$g$ & $9.8$ & \si{m/s^{2}} \\
$\mu_{g}$ & $1.8 \cdot 10^{-5}$ & \si{Pa.s} \\
$\mu_{l}$ & $1.516 \cdot 10^{-3}$ & \si{Pa.s} \\ 
$\epsilon$ & $10^{-8}$ & \si{m} \\
$L$ & $1000$ & \si{m} \\
\bottomrule
\end{tabular}
\end{table}


\subsubsection{Accuracy study}
\begin{figure}[hbtp]
\centering
	\begin{subfigure}[b]{.49\textwidth}
	\centering
	\includegraphics[width=\textwidth]{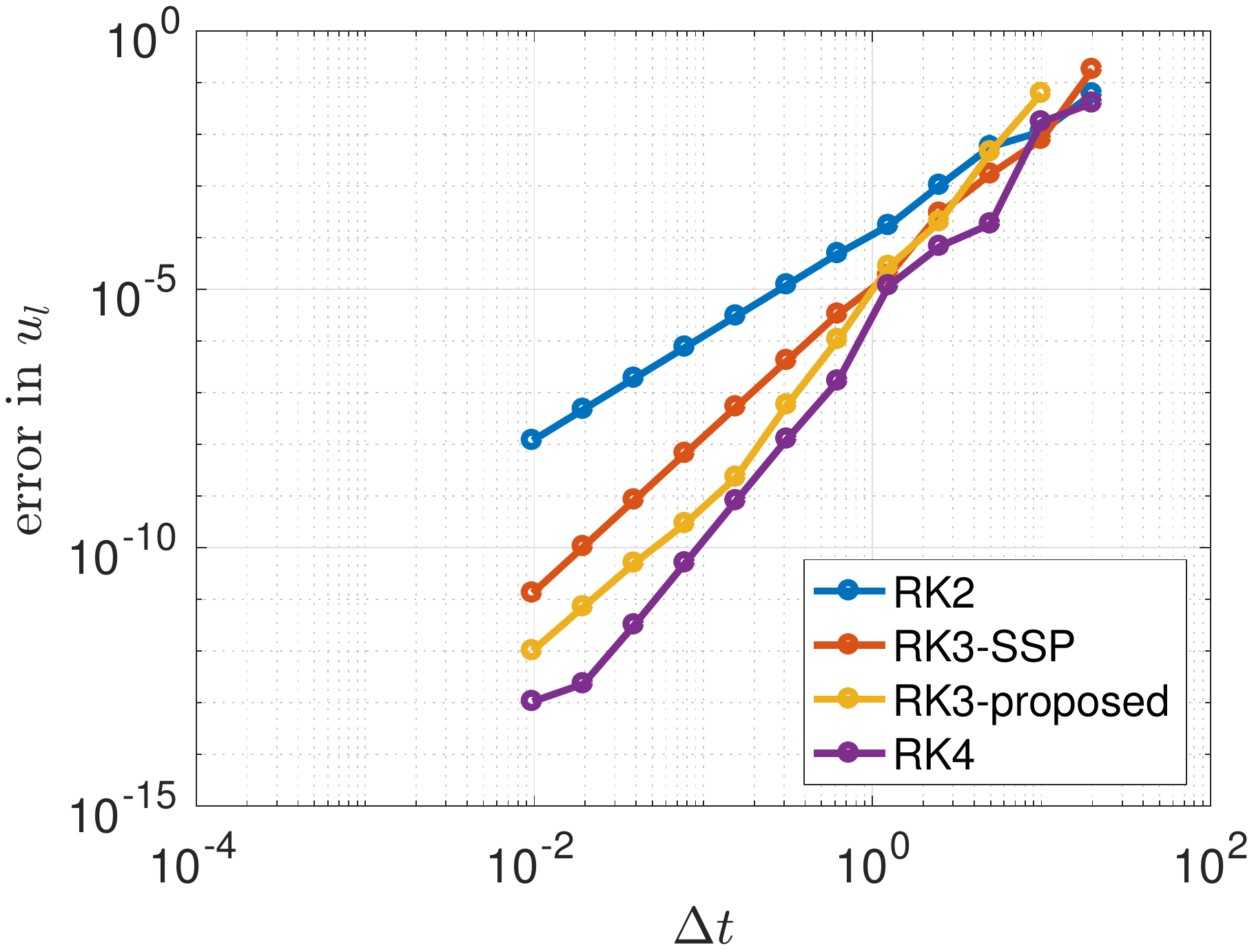}
	\caption{Weak boundary conditions.\label{fig:IFP_errors_a}}
	\end{subfigure}
	\hfill
	\begin{subfigure}[b]{.49\textwidth}	
	\centering
	\includegraphics[width=\textwidth]{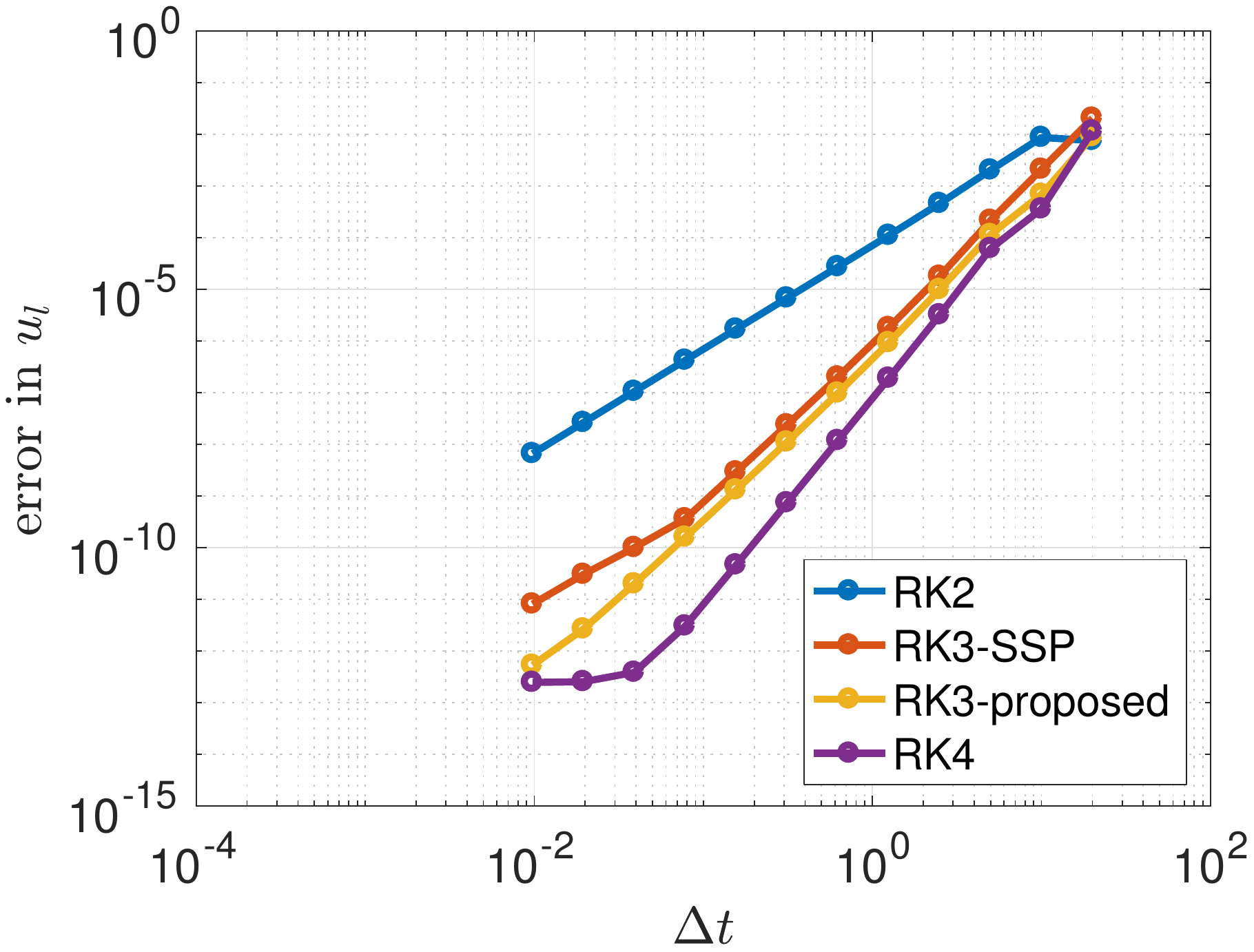}
	\caption{Strong boundary conditions.\label{fig:IFP_errors_b}}
	\end{subfigure}
\caption{Convergence of the temporal error in $u_{l}$ at $t=100$ \si{s} for the IFP problem.\label{fig:IFP_errors}}
\end{figure}

In order to show the accuracy of the time integration methods, we compare the solution at $t=100$ \si{s} for different time steps and Runge-Kutta methods in figure \ref{fig:IFP_errors}. The spatial grid is kept fixed at $N=40$ volumes, the time step varies from $\Delta t=20$ \si{s} to $\Delta t\approx 10^{-2}$ \si{s}, which means that the CFL number at the largest time step (based on the largest eigenvalue) is $\text{CFL} \approx 1.5$. A reference solution at $N=40$ with the HEM4 method \cite{Brasey1993} at very small time step ($\Delta t=10^{-3}$ \si{s}) is used to compute the temporal error, similar to equation \eqref{eqn:temporal_error}. 

Figure \ref{fig:IFP_errors} shows the error in $u_{l}$ for the case of weak (a) and strong (b) boundary imposition. For weak boundary conditions, all methods converge to the classic (non-DAE) order of convergence, since there are no additional order conditions. The irregular behaviour at coarse time steps might be attributed to the fact that the time step is of the same order as the period of the oscillation. Our proposed method, denoted RK3-proposed, converges with almost fourth order for coarse time steps, which can be attributed to the fact that the method satisfies most of the classical fourth order conditions. At small time steps, the third order behaviour is recovered. For strong boundary conditions, the effect of order reduction becomes apparent: the RK3-SSP method reduces to second order, whereas our proposed method does not suffer from order reduction. In this test case, the differential associated to the additional order condition (equation \eqref{eqn:additional_condition_3_differential}) is small, so the order reduction effect is only visible at small time steps. This is also the reason that RK4 does not show order reduction. In \ref{sec:MMS} we show that, for a different test case, also RK4 suffers from order reduction.


%
%
%

\FloatBarrier

\section{Conclusions}
A new constraint-consistent time integration strategy for the one-dimensional incompressible two-fluid model has been proposed. We have shown that the presence of the volume constraint in the two-fluid model manifests itself on the continuous, semi-discrete, and fully discrete level. On the continuous level, the volume constraint leads to Riemann invariants that correspond to hidden constraints of the model. On the semi-discrete level, the volume constraint leads to an index 3 differential-algebraic equation system, in which the same hidden constraints are present.

On the fully discrete level we have derived a novel time integration strategy, based on high-order `half-explicit' Runge-Kutta methods, that is consistent with these constraints. Our method is explicit for the mass and momentum equations and implicit for the pressure. The approach has a fractional-step like character: the pressure at a certain stage ($i$) of the Runge-Kutta method is such that the mixture velocity field is divergence free at the next stage ($i+1$), which in turn ensures that the phase masses satisfy the volume constraint at the following stage ($i+2$). Classic explicit Runge-Kutta methods can be used to achieve high order accuracy, provided that the boundary conditions are time-independent, or prescribed in a weak sense. For the important case of time-dependent boundary conditions prescribed in a strong manner, we have proposed a new three-stage, third order method which does not suffer from order reduction.

Our new time integration method has been demonstrated to perform according to the theoretical analysis for three problems, namely the Kelvin-Helmholtz instability for stratified pipe flow, sloshing in a closed pipe section, and ramp-up of the gas mass flow in a pipeline. It is shown that the classic fourth order Runge-Kutta method performs well for time-independent boundary conditions, but that our proposed third order method is the method of choice for time-dependent boundary conditions. 



Several important extensions of the current methodology are possible. 
First, a larger degree of implicitness can be necessary in certain problems. The current work has focused on half-explicit methods, in which only the pressure is computed implicitly. In case the two-fluid model is extended with terms involving short timescales, for example reaction terms due to thermodynamic phase transitions \cite{Lopez2017}, a higher degree of implicitness might be required. Possible extensions in line with the current work are IMEX (implicit-explicit) methods \cite{Ascher1997}, and partitioned or additive Runge-Kutta methods \cite{Kennedy2001}. 

Second, different type of constraint systems can be analyzed with our approach, such as the three-fluid model or the drift-flux model. This can shed new light on the wave structure of these models, see e.g.\ \cite{Evje2007}, and lead to improved time integration methods and boundary condition treatment. The compressible two-fluid model, see \cite{Sanderse2017}, does not possess the same constraint properties as the incompressible model, but the current analysis still provides an important limit that all-speed (incompressible - compressible) solvers should be able to handle.

Lastly, a great potential of our method lies in the extension to multi-dimensional problems, because the DAE analysis and proposed time integration strategy are still valid when the spatial dimension of the problem changes. Example application areas are multiphase flow in reactors \cite{Guelfi2007} and incompressible multiphase flow in reservoirs \cite{Bastian1999} (in which a saturation constraint similar to the volume constraint is present).




\section{Acknowledgements}
We would like to thank Ruud Henkes from Shell Technology Centre Amsterdam, who provided the initial inspiration to work on this topic, and valuable comments on the manuscript.

\appendix 

\section{Two-fluid model details}

\subsection{Geometry}\label{sec:geometry}
The following geometric identities are used to express the wall perimeters, interfacial perimeter, and liquid height in terms of the wetted angle $\gamma _l$ and pipe diameter $D=2 R$:
\begin{align}
{P_{gl}} &= D\sin {\gamma _l}, & {P_l} &= D{\gamma _l},\\
{P_g} &= D\left( {\pi  - {\gamma _l}} \right), & h &= \frac{1}{2} D\left( {1 - \cos {\gamma _l}} \right).
\end{align}
We use Biberg's approximation \cite{Biberg1999} to express $\alpha_l$ in terms of $\gamma_l$:

\begin{equation}\label{eqn:biberg}
\gamma_{l} = \pi \alpha_{l} + \left( \frac{3\pi }{2} \right)^{\frac{1}{3}} \left( \alpha_{g} - \alpha_{l} + \alpha _l^{\frac{1}{3}} - \alpha _g^{\frac{1}{3}} \right)   
                - \frac{1}{200} \alpha_{l} \alpha_{g} (\alpha_{g}-\alpha_{l}) (1 + 4(\alpha_{l}^2 + \alpha_{g}^2))
\end{equation}

\begin{figure}[h]
\fontfamily{lmss}\selectfont
\def\svgwidth{0.5\textwidth}
\centering
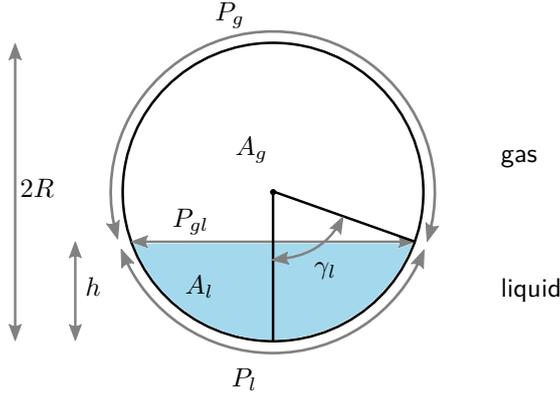
\caption{Stratified flow layout and definitions.}
\end{figure}

\subsection{Friction models}\label{sec:friction_models}
The wall (subscript $w$) and interfacial (subscript $gl$) shear stress are expressed by the Fanning friction factor definition:
\begin{equation}
\tau  = 
\begin{dcases} 
\frac{1}{2}{f_\beta }{\rho _\beta }{u_\beta }|{u_\beta }| & \text{wall} \\
\frac{1}{2}{f_{gl}}{\rho _g} ({u_\beta } - {u_\gamma })|{u_\beta } - {u_\gamma }| & \text{interfacial}
\end{dcases}
\end{equation}
The friction factor ${f_\beta }$ of phase $\beta$ with the pipe wall is modeled with the Churchill relation \cite{Churchill1977}:
\begin{align}
f_{\beta } &= 2{\left( {{{\left( {\frac{8}{{{{{\mathop{\rm Re}\nolimits} }_\beta }}}} \right)}^{12}} + {{\left( {A + B} \right)}^{ - 1.5}}} \right)^{\frac{1}{{12}}}}, \label{eq:churchill} \\
A &= {\left( {2.457\ln \left( {{{\left( {{{\left( {\frac{7}{{{{{\mathop{\rm Re}\nolimits} }_\beta }}}} \right)}^{0.9}} + 0.27\frac{\varepsilon }{{{D_{h\beta }}}}} \right)}^{ - 1}}} \right)} \right)^{16}},\\
B &= {\left( {\frac{{37530}}{{{{{\mathop{\rm Re}\nolimits} }_\beta }}}} \right)^{16}}.
\end{align}
Here $\varepsilon$ is the hydraulic pipe wall roughness, ${{\mathop{\rm Re}\nolimits} _\beta }$ is the Reynolds number,
\begin{equation}\label{eq:reynolds_number}
{{\mathop{\rm Re}\nolimits} _\beta } = \frac{{{\rho _\beta }{u_\beta }{D_{h\beta }}}}{{{\mu _\beta }}},
\end{equation}
and ${D_{h\beta }}$ is the hydraulic diameter:
\begin{equation}\label{eq:hydraulic_diameter}
{D_{h\beta }}  = 
\begin{dcases}
\frac{4{A_l}}{P_{l}} & \text{if }\beta  = l\\
\frac{4{A_g}}{P_{g} + P_{gl}} & \text{if }\beta  = g
\end{dcases}
\end{equation}
The interfacial friction factor ${f_{gl}}$ is calculated by \cite{Liao2008}:
\begin{equation}\label{eq:interface_friction}
f_{gl} = \max (f_g,0.014).
\end{equation}
\section{Order reduction analysed with the method of manufactured solutions}\label{sec:MMS}
In the method of manufactured solutions (MMS) an analytical solution is assumed, substituted into the two-fluid model equations, and the resulting term is used as an additional known source term in the two-fluid model. We design an analytical solution $W^{*}(s,t)$ which can be exactly represented by the spatial discretization, so that any errors in the numerical solution are purely due to the time integration method.
For this purpose, the phase masses $m_{g}$, $m_{l}$ are chosen to be constant in space, and time-varying according to a prescribed function $f(t)$:
\begin{align}
m_{g}^{*} (s,t) &= m_{g}^{*} (t) = \rho_{g} A_{g}^{*} (t) =  \rho_{g} \hat{A}_{g} f(t),\\
m_{l} ^{*}(s,t) &= m_{l}^{*} (t) = \rho_{l} A_{l}^{*} (t) =  \rho_{l} (A - \hat{A}_{g} f(t)). 
\end{align}
The phase momenta are chosen to be linearly varying in space, in such a way that no source term appears in the mass equations:
\begin{align}
I_{g}^{*}(s,t) &= m_{g}^{*}(t) \hat{u}_{g} - \dot{m}_{g}^{*} (t) s = \rho_{g} \hat{A}_{g}  \left( \hat{u}_{g} f(t) - \dot{f}(t) s \right), \\
I_{l}^{*}(s,t) &=  m_{l}^{*}(t) \hat{u}_{l} - \dot{m}_{l}^{*} (t) s = \rho_{l} (A_{l}^{*}(t) \hat{u}_{l} + \hat{A}_{g} \dot{f}(t) s). 
\end{align}
The phase velocities follow from $u_{g}^{*}(s,t)=I_{g}^{*}(s,t)/m_{g}^{*}(t)$ and similarly for $u_{l}^{*}(s,t)$. Note that $u_{g}^{*}(s=0,t)$ is constant in time, although $I_{g}^{*}(s=0,t)$ is not. For the pressure a linearly varying profile in space is assumed,
\begin{equation}
p^{*}(s) = c_{1} s + c_{2}.
\end{equation}
Time dependency could be incorporated in the pressure solution, but is not very important due to the Lagrange-multiplier nature of the pressure.


The additional source term $F_{\text{body}}$ required in the momentum equation that forces this analytical solution is then given by
\begin{equation}
\begin{split}
F_{\text{body},g} (s,t) &= \frac{\partial }{\partial t} \left( I_{g}^{*}(s,t) \right) + \frac{\partial }{\partial s} (I_{g}^{*}(s,t) u_{g}^{*}(s,t))  + \frac{\partial p^{*}(s) }{\partial s} A_{g}^{*}(t) + S_{g}^{*}(s,t),\\
&=  \rho_{g} \hat{A}_{g} (\hat{u}_{g} \dot{f}(t) - \ddot{f}(t) s) + 2 \rho_{g} \hat{A}_{g} \left( s \frac{\dot{f}^{2}(t)}{f(t)}  - \hat{u}_{g} \dot{f}(t) \right) + \hat{A}_{g} f(t) c_{1} + S_{g}^{*}(s,t),
\end{split}
\end{equation}
where
\begin{equation}
S_{g}^{*} (s,t) = S_{g}( u_{g}^{*}(s,t), u_{l}^{*}(s,t),A_{g}^{*}(t))
\end{equation}
contains the algebraic source terms (friction and gravity). The level gradient term is zero because the hold-up fractions are uniform in space. The source term for the liquid momentum equation is constructed in a similar fashion.

The function $f(t)$ is chosen to be continuously differentiable with non-vanishing derivatives in order to prevent false perception of high order accuracy. We therefore choose the function
\begin{equation}\label{eqn:MMS_f}
f(t) = \frac{1}{60} (\sin(a t)+5) e^{b t}, \quad a=2, \quad b=1/20.
\end{equation}
The amplitudes $\hat{A}_{g}$ and $\hat{u}_{g}$ are chosen according to the steady state solution.

\begin{figure}[hbtp]
\centering
	\begin{subfigure}[b]{.49\textwidth}
	\centering
	\includegraphics[width=\textwidth]{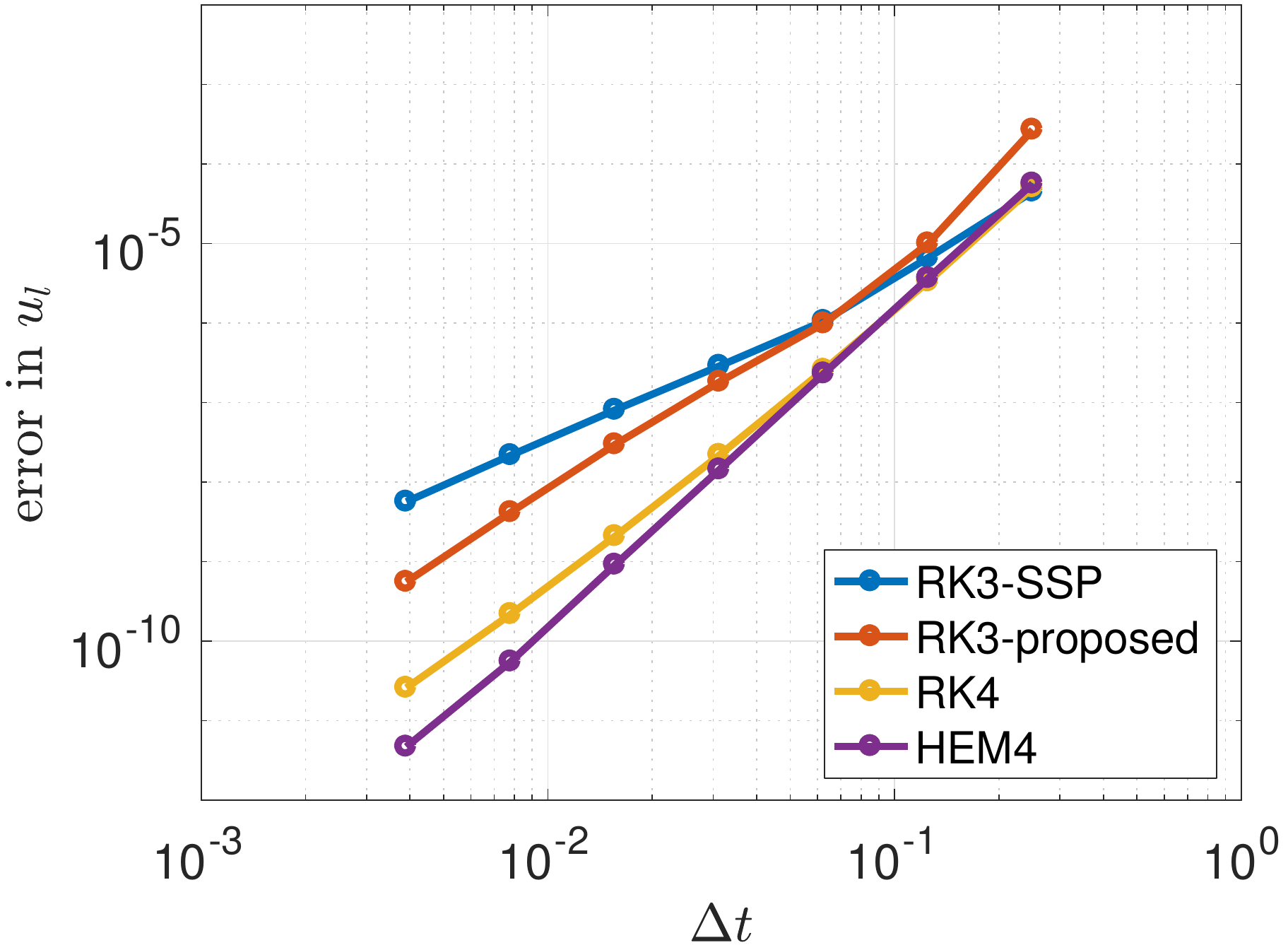}
	\caption{$u_l$, strong BC.}
	\end{subfigure}
	\hfill
	\begin{subfigure}[b]{.49\textwidth}	
	\centering
	\includegraphics[width=\textwidth]{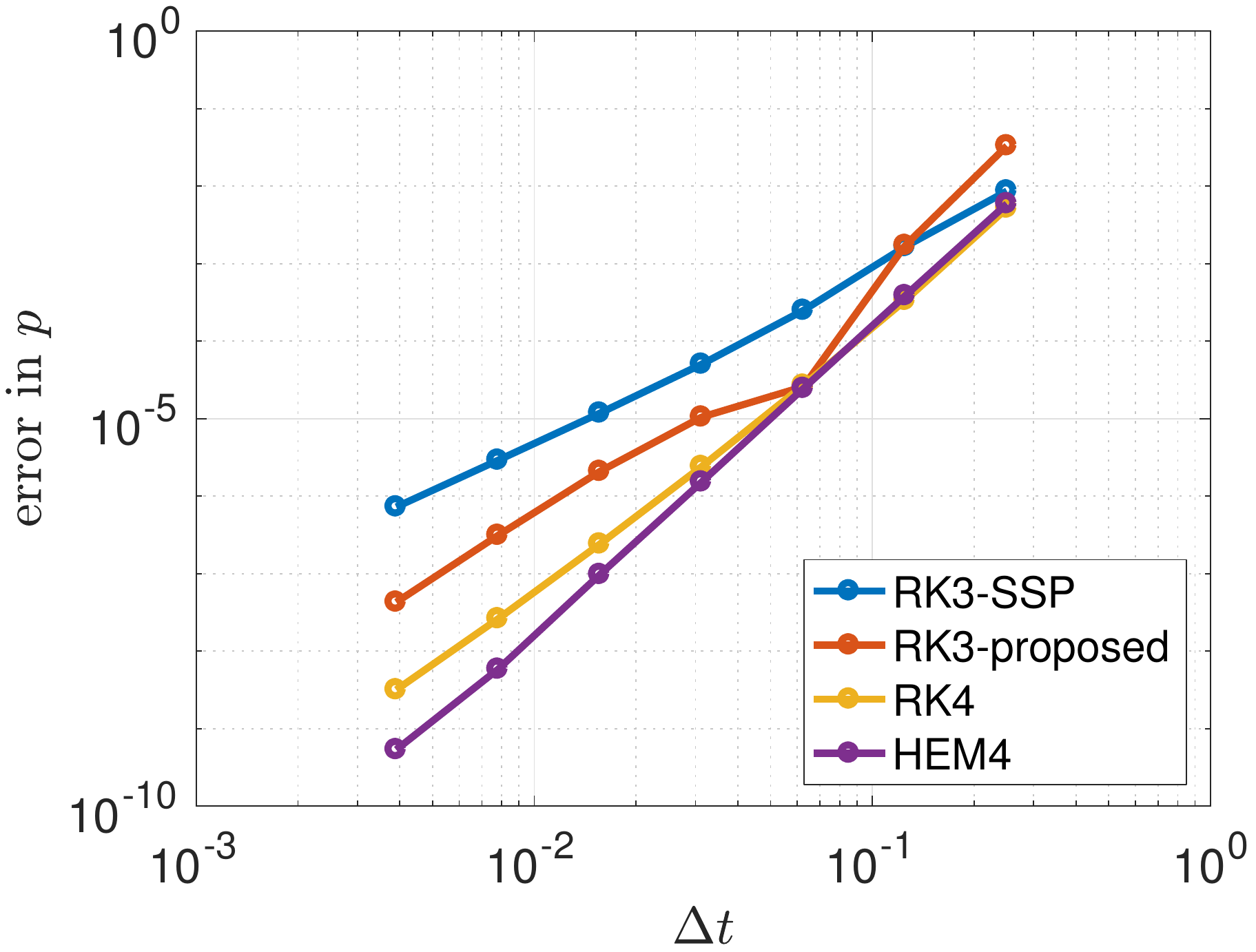}
	\caption{$p$, strong BC.}
	\end{subfigure}\\
	\vspace{0.5cm}
	\begin{subfigure}[b]{.49\textwidth}
	\centering
	\includegraphics[width=\textwidth]{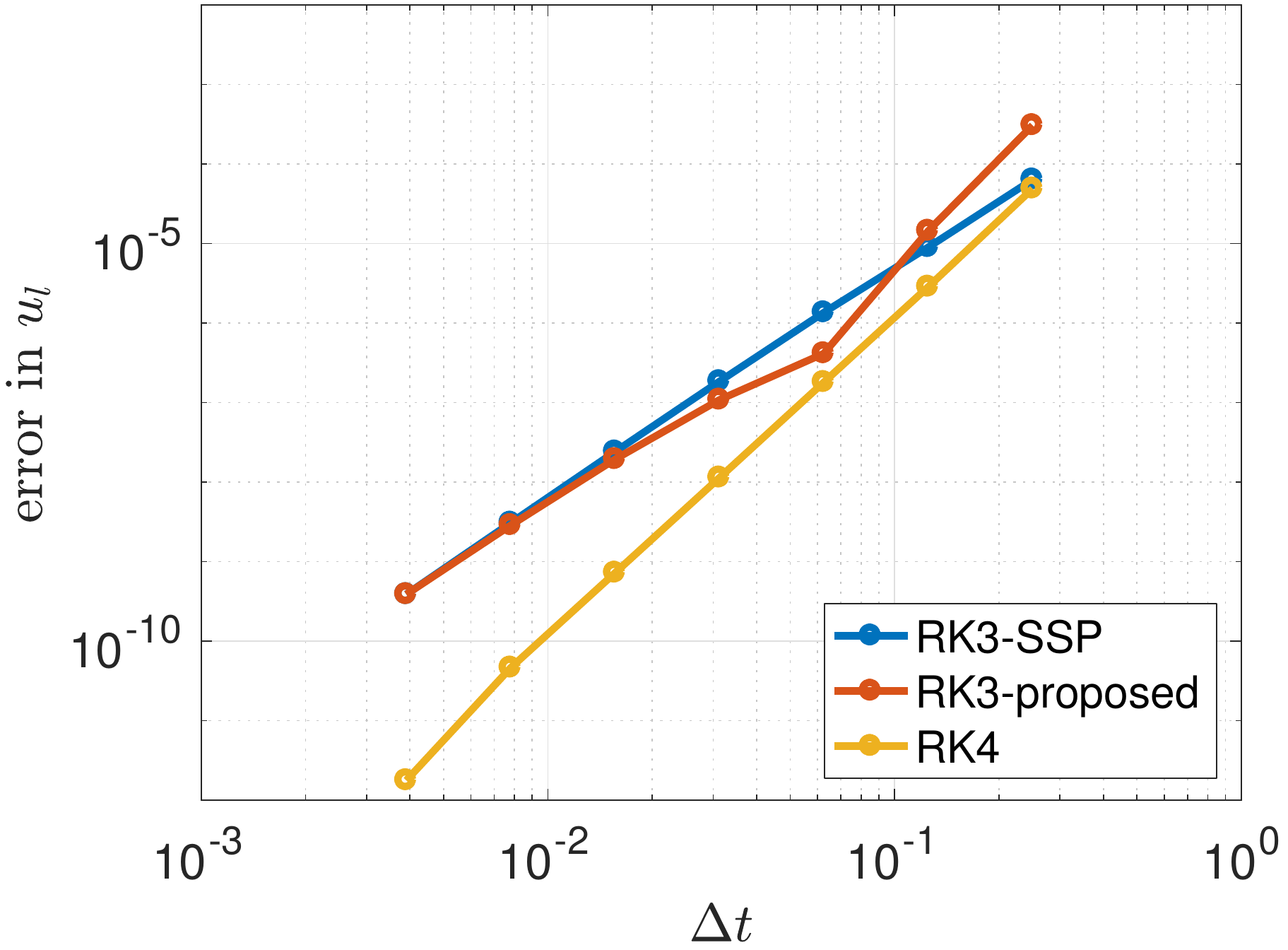}
	\caption{$u_l$, weak BC.}
	\end{subfigure}
	\hfill
	\begin{subfigure}[b]{.49\textwidth}	
	\centering
	\includegraphics[width=\textwidth]{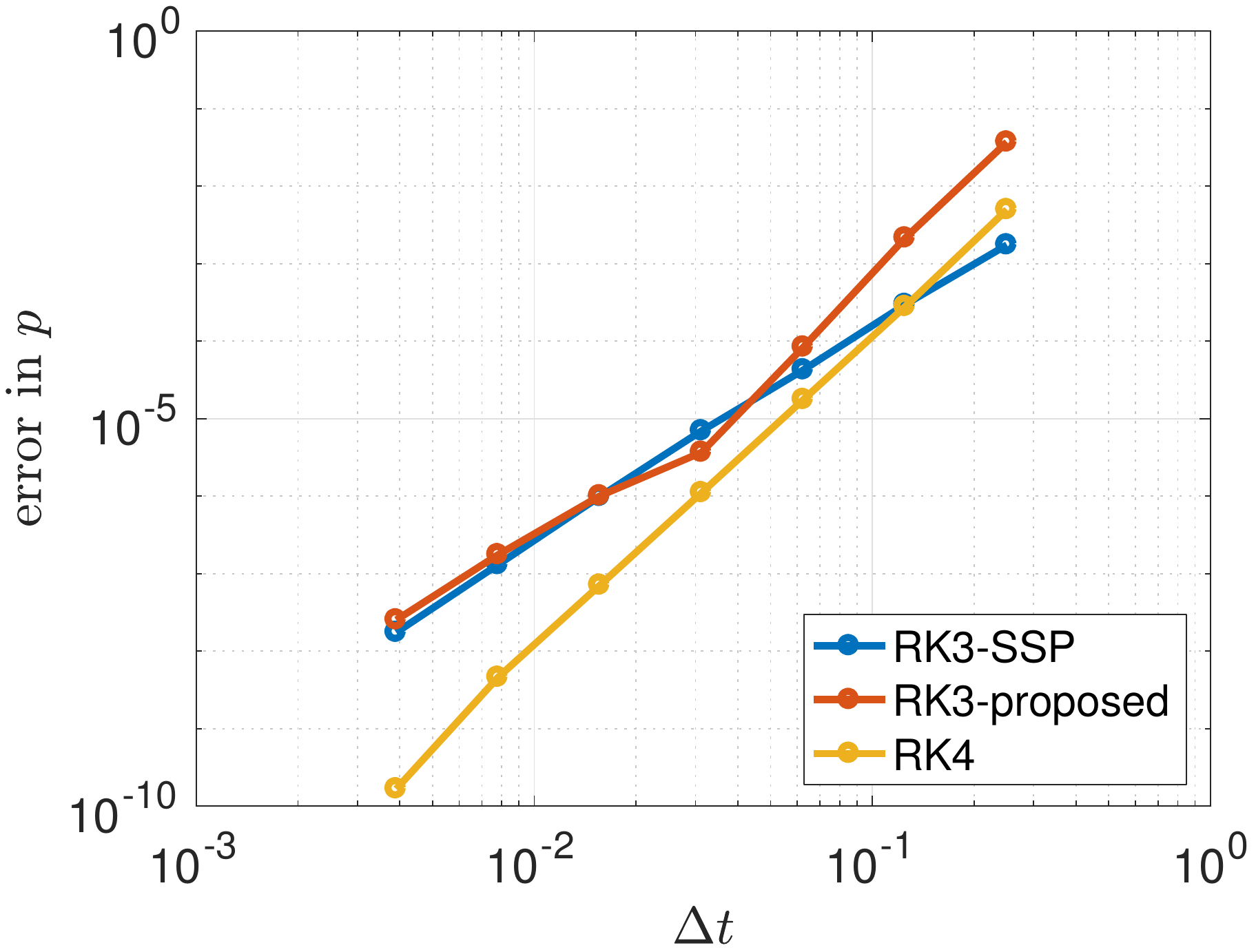}
	\caption{$p$, weak BC.}
	\end{subfigure}
\caption{Convergence of the temporal error for MMS. \label{fig:MMS_errors}}
\end{figure}

The parameters of the test case are the same as in section \ref{sec:IFP}, except that the pipe length is $L=10$ \si{m}, the diameter is $D=0.25$ \si{m},  laminar friction closure is used, and  the initial gas and liquid flow rates are $0.04$ and $2$ \si{kg/s}, respectively. At $s=0$ unsteady Dirichlet conditions according to \eqref{eqn:MMS_f} are prescribed; at $s=L$ outflow conditions are used. We integrate the two-fluid model equations until $t=20$ \si{s}. Figure \ref{fig:MMS_errors} shows the temporal errors for the liquid velocity and the pressure, for both strong and weak boundary conditions. In this test case it is evident that the classic RK schemes RK3-SSP and RK4, not designed for DAEs, suffer from order reduction: RK3-SSP reduces to second order, and RK4 to third order. Our proposed RK3 scheme remains third order. The HEM4 method of \cite{Brasey1993} is fourth order, but requires five stages, and is therefore less attractive.

\bibliography{refs}
\bibliographystyle{abbrv}

\end{document}